\documentclass[11pt]{amsart}

\usepackage{theoremref}
\usepackage{mathrsfs}
\usepackage{hyperref}

\usepackage{tikz}
\usepackage{times}
\usepackage{color}
\usepackage{amsmath}
\usepackage{amssymb}
\usepackage{amsthm}
\usepackage{enumerate}
\usepackage{amsbsy}
\usepackage{amsfonts}
\newtheorem{theo}{Theorem}[section]
\newtheorem{defin}[theo]{Definition}
\newtheorem{prop}[theo]{Proposition}

\newtheorem{lemm}[theo]{Lemma}
\newtheorem{rem}[theo]{Remark}

\newcommand{\al}{\alpha}
\newcommand{\be}{\beta}
\newcommand{\ga}{\gamma}
\newcommand{\Ga}{\Gamma}

\newcommand{\ep}{\epsilon }
\newcommand{\te}{\theta}
\newcommand{\De}{\Delta}
\newcommand{\de}{\delta}

\newcommand{\pa}{\partial}

\newcommand{\R}{{\mathbb R}^n}
\newcommand{\hR}{{\mathbb R}^n_+}

\newcommand{\ri}{\rightarrow}
\newcommand{\Rn}{{\mathbb R}^{n-1}}

\newcommand{\na}{\nabla}

\newcommand{\calF}{{ \mathcal F }}
\newcommand{\calG}{{ \mathcal G }}
\newcommand{\calR}{{ \mathcal R }}

\newcommand{\bke}[1]{\left( #1 \right)}
\newcommand{\bkt}[1]{\left[ #1 \right]}
\newcommand{\bket}[1]{\left\{ #1 \right\}}
\newcommand{\norm}[1]{\left\Vert #1 \right\Vert}
\newcommand{\abs}[1]{\left| #1 \right|}
\newcommand{\rabs}[1]{\left. #1 \right|}

\newcommand{\calB}{{ \mathcal B  }}

\newcommand{\hoe}{H\"{o}lder }

\newtheorem{assumption}{Assumption}
\newtheorem{notation}{Notation}

\begin{document}
\baselineskip=18pt

\title[Unbounded normal derivatives]{
 Singular weak solutions near boundaries  in a half space away from localized force for the Stokes and Navier-Stokes equations}
\author{Tongkeun Chang and Kyungkeun Kang}

\thanks{}

\begin{abstract}
We prove that there exists a weak solution of the Stokes system with a non-zero external force and no-slip boundary conditions in a half space of dimensions three and higher so that its normal derivatives are unbounded near boundary. A localized and divergence free singular force causes, via non-local effect, singular behaviors of normal derivatives for the solution near boundary, although such boundary is away from the support of the external force.
The constructed one is a weak solution that has finite energy globally, and it can be comparable to the one in \cite{Seregin-Sverak10} as a form of a shear flow that is of only locally finite energy. Similar construction is performed for the Navier-Stokes equations as well.


\end{abstract}

\maketitle

\section{Introduction}
\setcounter{equation}{0}

In this paper, we consider the following non-stationary Stokes system with non-zero external force, zero initial data and no-slip boundary condition in half space $\R_+$, $n \ge 3$;
\begin{align}\label{maineq}
\left\{\begin{array}{l} \vspace{2mm}
w_t -\De w + \na \Pi  = f,\quad  {\rm div} \, w =0 \qquad \R_+ \times (0, 1),\\
\vspace{2mm}
\rabs{w}_{x_n =0} =0, \qquad \rabs{w}_{t =0} =0.
\end{array}
\right.
\end{align}
Here we assume that $f$ is compactly supported in $\overline{\R_+} \times (0, 1)$.
Our concern is local analysis of the solution of the Stokes system \eqref{maineq}
near boundary, in particular, in the region near boundary away from the support of $f$.
A specific form of localized external force $f$ in \eqref{maineq} is described in {\bf Assumption \ref{force-f}}.
One can imagine similar situation for the heat equation in a half space
\begin{align*}
\left\{\begin{array}{l} \vspace{2mm}
u_t -\De u = f, \qquad \R_+ \times (0, 1),\\
\vspace{2mm}
\rabs{u}_{x_n =0} =0, \quad \rabs{u}_{t =0} =0.
\end{array}
\right.
\end{align*}
We suppose $f$ is compactly supported, e.g., in $\overline{B^+_1}\times (0,1)$, where $B^+_r=\bket{x\in \mathbb R^{n}: |x|<r, x_n>0}$. Even in the case that $f$ is singular in $B^+_1\times (0,1)$, it is known due to classical regularity theory that $u$ becomes regular, in particular, near the boundary, away from the support of $f$, namely
\begin{equation}\label{localheatineq}
\norm{\partial_t^m \partial_x^l u}_{L^{\infty}(B^+_{x',r}\times (t-r^2, t))} \le c\norm{u}_{L^{2}(B^+_{x',2r}\times (t-4r^2, t))},\qquad m, l \ge 0,
\end{equation}
where $x' \in  \bket{y'\in \mathbb R^{n-1}: |y'|>2}$ and $ (t-4r^2, t)\subset (0,1)$.
It is, however, unclear, due to the nonlocal effect, whether or not such an estimate is available for the Stokes system \eqref{maineq}.

One can compare to the
Stokes system with nonzero boundary data, instead nonzero force, in half space $\R_+$, that is
\begin{align}\label{SS-bdata}
\left\{\begin{array}{l} \vspace{2mm}
w_t -\De w + \na \Pi  = 0,\quad  {\rm div} \, w =0, \qquad \R_+ \times (0, 1),\\
\vspace{2mm}
\rabs{w}_{x_n =0} =\varphi (x', t), \quad \rabs{w}_{t =0} =0.
\end{array}
\right.
\end{align}
In this case, for the localized boundary data, it has been shown that the estimate
\eqref{localheatineq}  is, in general, not true for the Stokes system \eqref{SS-bdata}, and
 furthermore, construction of solutions with the same singular behaviors has been also constructed for the Navier-Stokes equations as well (see \cite{CKca}, \cite{CK0206}, \cite{Kang05}, \cite{KLLT21}).
In particulr, it was shown that the constructed singular solutions in \cite{KLLT21} are indeed global energy solutions, i.e. $w\in L^{\infty}((0,1);L^2 (\R_+))\cap L^2((0,1); \dot{H}(\R_+))$, $n \ge 3$, for the Stokes system and Navier-Stokes equations as well.
Therefore, we can say that,
unlike the heat equation,  non-local effect of the Stokes system with singular non-zero boundary data may cause  violation of local smoothing effects of solutions near the boundary.
However, most of examples have been constructed via nonzero flux at the boundary and it is not clear whether or not singular behaviors of solutions with finite global energy can be developed in the case of no-slip boundary condition on the every boundary.
Nevertheless, Seregin and \u{S}ver\'ak found in \cite{Seregin-Sverak10} the form of shear flow whose
 normal derivatives of solutions are unbounded near boundary in the half-space. More precisely, in \cite{Seregin-Sverak10}, they constructed the following form of shear flow:
\[
w(x,t)=(u(x_3,t), 0,0),\quad \,\,\Pi(x,t)=-g(t)x_1,\quad \mathbb R^3_+\times (-4,0)
\]
with homogeneous initial and boundary conditions and $g(t)=\abs{t}^{-1+\alpha}$, $\alpha\in (0, \frac{1}{2})$.
Then the solution is explicitly given as
\begin{equation}\label{SS-example}
w(x_3,t)=\frac{2}{\sqrt{\pi}}\int_{-4}^t g(t-\tau-4)d\tau\int_0^{\frac{x_3}{\sqrt{4(\tau+4)}}} e^{-\xi^2}d\xi
\end{equation}
and one can see that $w$ is bounded but
$\partial_{x_3} w\ge Cx^{-1+2\alpha}_3$ in the region near $x_3=0$ with $x^2_3>-4t$. We remark that the solution is not of finite energy in the half space and it is not even decaying, as $x_3$ tends to infinity.

 Our main objective of this paper is to construct solutions of \eqref{maineq} such that it is of finite energy in the half space, i.e. global weak solutions, and however it has a singular behavior near boundary,  namely, unbounded normal derivatives such that main features of blow-up profiles are similar to the one \eqref{SS-example} specified in \cite{Seregin-Sverak10}.

Firstly, we specify the external force $f$, which is of divergence free and has a certain type of singular behaviors for normal variable $x_n$ and time variable $t$.
For convenience, we denote $x=(x', x_n )\in \mathbb R^n$ with $x'\in \mathbb R^{n-1}$.

\begin{assumption}\label{force-f}
Let $n\ge 3$ and $0 <\al, \, \be < 1$. Suppose that $g:\mathbb R^n_+ \rightarrow \mathbb R$ is a real-valued function of the form $g (x) = g^{{\mathcal T}}(x') g^{{\mathcal N}}(x_n)$, where non-negative functions $ g^{{\mathcal T}}: \mathbb R^{n-1} \rightarrow \mathbb R$ and $g^{{\mathcal N}}: \mathbb R_+ \rightarrow \mathbb R$ satisfy
\begin{equation*}
g^{{\mathcal T}} \in C_c^\infty(\mathbb R^{n-1}), \qquad
{\rm supp}\, g^{{\mathcal T}} \Subset B^{'}_1=\bket{x'\in\mathbb R^{n-1}: |x'|<1},
\end{equation*}
\begin{equation*}
g^{{\mathcal N}} \in C^\infty(\mathbb R_+), \quad {\rm supp} \, g^{{\mathcal N}} \subset  (0, 2), \quad g^{{\mathcal N}}(x_n) = x_n^{1 -\be}, \,\,\mbox{ for }\,\, x_n \in (0, 1].
\end{equation*}
Let $a>0$ be a constant.
Furthermore, we suppose that a vector field $f=(f_1, \cdots, f_n):\mathbb R^n_+ \times \mathbb R_+ \rightarrow \mathbb R^n$ is given as $f_2=a\frac{\pa g}{\pa x_n} (x) h(t)$, $f_n=-a\frac{\pa g}{\pa x_2} (x)  h(t)$ and $f_i=0$ for $i\neq 2$, $i\neq n$, i.e.
\begin{equation}\label{force-st}
f (x,t) = \bke{0,  a\frac{\pa g}{\pa x_n} (x) h(t),0, \cdots, 0, -a\frac{\pa g}{\pa x_2} (x)  h (t)},
\end{equation}
where a non-negative function $h:\mathbb R_+ \rightarrow \mathbb R$ is given by
\begin{align*}
 h(t) = (t-\frac12)^{-\al} \chi_{(\frac12,\infty)}(t).
\end{align*}
\end{assumption}

\begin{rem}
We note that the vector field $f$ in   \eqref{force-st} is divergence free in $\mathbb R^n_+$ and the normal component vanishes at the boundary, namely
${\rm div} \, f =0$ and $\rabs{f_n}_{x_n =0} =0$. It is straightforward that
\begin{equation}\label{mixednorm-f}
f \in L^{q_1}_t L^{p_1}_x (\R_+ \times (0, \infty)),\qquad q_1 \in [1, \frac{1}{\al}),\quad p_1 \in [1, \frac{1}{\be}).
\end{equation}
We remind that $f_2$ and $f_n$ near $(x_n, t)=(0, 1/2)$ behave as follows:
\begin{equation*}
f_2 \sim x^{-\beta}_n \bke{t-\frac{1}{2}}^{-\alpha},\qquad
f_n \sim x^{1-\beta}_n \bke{t-\frac{1}{2}}^{-\alpha}.
\end{equation*}
\end{rem}

It was shown in \cite{So1} that, in case that a given vector field $f$ in a half space holds ${\rm div} \, f =0$ and $\rabs{f_n}_{x_n =0} =0$,
the solution  $w$ and the associated pressure $\Pi$ in \eqref{maineq} are represented by
\begin{align}\label{int-0909-1}
w(x,t) = \int_0^t \int_{\R_+} K(x,y, t-s)  f(y,s) dyds,
\end{align}
\begin{equation}\label{int-expression-p}
\Pi(x,t) =\int_0^t \int_{{\mathbb R}^n_+} P(x,y, t-s)
\cdot f(y,s)dyds,
\end{equation}
where the Green tensor $K=(K_{ij})$ and the pressure vector $P=(P_j)$ are given as
\begin{align}\label{formulas-k}
\notag K_{ij}(x,y, t) &=\de_{ij} \big(\Ga(x-y,t) - \Ga(x -y^*,t) \big)\\
& \qquad  - 4(1 -\de_{jn}) D_{x_j} \int_0^{x_n} \int_{\Rn}\Ga(x -y^* -z, t) D_{z_i} N(z) dz,
\end{align}
\begin{align}\label{formulas-p}
P_j(x,y,t) =4 (1 - \de_{jn}) D_{x_j}(D_{x_n}+D_{y_n})\int_{{\mathbb
R}^{n-1}} N(x- z')\Ga(z' -y, t) dz',
\end{align}
where $y^* = (y', -y_n)$, $\de_{jn}$ is Kronecker delta function, and $\Gamma(x,t)=(4\pi t)^{-\frac{n}{2}}e^{-\frac{|x|^2}{4t}}$ and $N(x)=-c_n|x|^{2-n}$ with $c_n =(n(n-2)\omega_n)^{-1}$ denote Gaussian kernel and Newtonian kernel in $n$ dimensions, $n \ge 3$, respectively.

Next, we introduce the notion of weak solution for the Stokes system \eqref{maineq}.

\begin{defin}\label{weaksol-SS}
Let $T\in (0, \infty)$ and $ f \in L^q_tL^p_x(\mathbb R^n_+ \times (0,T))$ for $1< p,q<\infty$. We say that a vector field $w\in L^2(0,T; \dot H^1 (\R_+)) \cap L^\infty(0,T; L^2 (\R_+))$ is a weak solution of the Stokes system \eqref{maineq},
if the following equality is  satisfied:
\begin{align}\label{weakform-SS}
\int^T_0\int_{\hR }\nabla w : \nabla \Phi
dxdt=\int^T_0\int_{\hR}\bke{w \cdot \Phi_t+ f \cdot \Phi} dxdt
\end{align}
for every vector field $\Phi\in C^2_c( \hR\times [0, T))$ with $\mbox{\rm div } \Phi=0$, and
in addition, for every scalar function $\Psi\in C^1_c(\overline{\hR})$
\begin{equation*}
\int_{\hR} w(x,t) \cdot \na \Psi(x) dx =0  \quad  \mbox{ for all}
\quad 0 < t< T.
\end{equation*}
Furthermore, for every vector field $\varphi\in C^0_c(\overline{\hR})$
\begin{equation*}
\lim_{t\rightarrow 0}\int_{\hR} w(x,t) \cdot\varphi(x) dx=0
\end{equation*}
\end{defin}
From now on, for simplicity, we assume that $T=1$,  without loss of generality, in Definition \ref{weaksol-SS}.
The concept of weak solutions can be relaxed by removing the restriction that solutions belong to energy class, i.e. $  L^2(0,1; \dot H^1 (\R_+)) \cap L^\infty(0,1; L^2 (\R_+))$. Indeed, for comparison, we also introduce a notion of {\it very weak solutions} (see Definition \ref{stokesdefinition} in Section 2).

Our main objective of the paper is to construct a weak solution of the Stokes system \eqref{maineq}
with singular behavior near boundary. To be more precise, normal derivativies of weak solutions are unbounded at the boundary away from the support of $f$, although solutions are in energy classes and even locally bounded.

\begin{notation}\label{boundaryset}
Let $i$ be an integer with $1\le i\le n-1$ and $i\neq 2$.
We introduce, for convenience, a set $A_{i}\subset \Rn$ defined by
\begin{equation}\label{Ai-set}
A_{i}=\bket{x'\in \Rn:
\frac12 | x_i| \leq |x_2| \leq 2|x_i|,\,\,
|x'|^2 \leq 2 \big( |x_i|^2 + |x_2|^2), \,\, |x_i|, |x_2|>2}.
\end{equation}
We split $A_i$ into two disjoint sets, denoted by $A_{i1}$ and $A_{i2}$, as follows:
\begin{equation}\label{Ai12-sets}
A_{i1}=A_i \cap\bket{x'\in \Rn:  x_i x_2 > 0}, \qquad
A_{i2}=A_i \cap\bket{x'\in \Rn:  x_i x_2 < 0}.
\end{equation}
We also denote $B_i:=B_{i1} \cup B_{i2}$, where $B_{i1}$ and $B_{i2}$ are defined by
\begin{equation}\label{Bi1-set}
B_{i1} = \bket{x' \in \Rn \,\big| \,  \frac{1}{4\sqrt{n} }\abs{x'} > |x_2|, \,\,  2<|x_i| < \infty   },
\end{equation}
\begin{equation}\label{Bi2-set}
B_{i2} = \bket{x' \in \Rn \,\big| \,  4\sqrt{n}   \abs{x'} <  |x_2|, \,\, 2< |x_2| < \infty   }.
\end{equation}
\end{notation}

A two dimensional cartoon of sets defined above, for example, is pictured in the Appendix \ref{boundaryset-1}.

Now we are ready to state first main result.

\begin{theo}\label{maintheo0503}
Let   $f$ be given in Assumption \ref{force-f}, and $A_i$ and $B_i$ disjoint sets
defined in Notation \ref{boundaryset}.
\begin{itemize}
\item[(i)]
If $0 <\be<\frac12$, then the solution $w$ defined in \eqref{int-0909-1} of the Stokes system \eqref{maineq} becomes the weak solution satisfying
\begin{equation}\label{energy-bound}
\| w\|_{L^\infty_t L^2_x(\R_+ \times (0, 1))} +\|D_{x} w\|_{L^2 (\R_+ \times (0, 1))}  \le c=c(
\norm{f}_{L^{q_1}_t L^{p_1}_x (\R_+ \times (0, 1))}),
\end{equation}
where $q_1 \in [1, \frac{1}{\al})$ and $p_1 \in [1, \frac{1}{\be})$.
\item[(ii)]
If $ q > 6$ and  $2 + \frac3{q}   < 2\al +\be$, then  normal derivatives of $w$ are  singular on any subset of $A_i \cup B_i$, i.e.
\begin{equation}\label{0504-2}
\|\partial_{x_n} w\|_{L^l  (D\times (0,1)\times (0,1) )}
= \infty,  \qquad \mbox{ for any } l \geq q \mbox{ and  } D\subset A_i \cup B_i ,
\end{equation}
\end{itemize}
\end{theo}

\begin{rem}
\begin{itemize}
\item[(i)]
We note that in case  $ 6 < q < \infty$,  there are $\al\in (0,1)$ and $\beta\in (0, 1/2)$ such that
$2 + \frac3{q}   < 2\al +\be$, and thus
\begin{align*}
 \bket{ (\al, \be) \in (0, 1) \times (0, \frac{1}{2}) \,: \, 2 + \frac3{q}   < 2\al +\be} \neq \emptyset.
\end{align*}
\item[(ii)]
We remark that if we don't require that solutions belongs to the energy class and
instead, if we allow it to be a very weak solution (see Definition \ref{stokesdefinition}), then
it is unnecessary to assume that $0 <\be<\frac12$, and thus it is possible to construct a very weak solution $u$ such that $\nabla u$ becomes unbounded in $L^q_{\rm{loc}}$, $q>3$. near boundary away from the support of $f$.
It turns out that such examples show similar singular behaviors as those of the example constructed in \cite{Seregin-Sverak10}. Since our concern is about weak solutions, we are not going to pursue the matter on very weak solutions in this paper.
\end{itemize}
\end{rem}

Secondly, we similarly analyze pressure both globally and locally, and
obtain the following:

\begin{theo}\label{maintheo0503-2}
Let $\al$, $\be$ and $f$ be the numbers and the vector field  in Theorem \ref{maintheo0503}.
Let $\Pi$ be a  pressure associated with the weak solution $w$ of the Stokes system in Theorem \ref{maintheo0503}, defined by \eqref{int-expression-p}.
\begin{itemize}
\item[(i)]
Let $p>\frac{2n}{n-1}$ and $q>1$. Suppose that $\alpha \in (0,1)$ and $\be\in (0, \frac12)$  are numbers satisfying
\begin{align*}
\be > \frac{n}{(n-1)p},\qquad 2\al +n \be<\frac2q  +\frac{n}{p} + 1.
\end{align*}
Then, the pressure $\Pi$ is globally bounded in $L^{q}_t L^p_x$ and satisfies
\begin{equation}\label{pi-bound}
\| \Pi\|_{L^{q}(0, 1; L^p (\R_+))} \le c=c(
\norm{f}_{L^{q_1}_t L^{p_1}_x (\R_+ \times (0, 1))}),
\end{equation}
where $q_1 \in [1, \frac{1}{\al})$ and $p_1 \in [1, \frac{1}{\be})$.
\item[(ii)]
If $q\in (1, \infty)$ satisfies
\begin{align}\label{0812-2}
1 +\frac2q < 2\al +\be,
\end{align}
then $\Pi$ is locally unbounded in $L^{q}_{x,t}$, that is
\begin{equation}\label{0504-4}
\| \Pi\|_{L^{q} ( \{ |x'| > 2 \} \times (a,b)\times (0, 1))} =\infty %\qquad 0 \leq a < b < \infty.
\end{equation}
for any $a, b$ with $0 \leq a < b < \infty$.
\end{itemize}
\end{theo}

\begin{rem}
As mentioned early, in the proof of Theorem \ref{maintheo0503},   $\beta\in (0, \frac{1}{2})$ and $ \al \in (0, 1)$ are imposed, which implies $ q> \frac43$ in (ii) of Theorem \ref{maintheo0503-2}.
If it is, however, extend to the very weak solutions, the condition $\beta\in (0, \frac{1}{2})$ can be relaxed as $\beta\in (0, 1)$, and thus \eqref{0812-2} is valid for any $q>1$. In addition,
the condition $p>\frac{2n}{n-1}$ in (i)  can be relaxed by $p>\frac{n}{n-1}$ for such case.
\end{rem}

\begin{rem}
It is not difficult to see that \eqref{pi-bound} and  \eqref{0504-4} are compatible. Indeed,
suppose that $ \frac{2n}{n-1} < p < \infty $, $1 < q < \infty$. Then, we can check easily that two sets $C$ and $D$ below have no intersection, i.e. $C \cap D  = \emptyset$, where
\begin{align*}
C & = \bket{ (\al, \be) \in (0, 1) \times (0, 1) \,\,:\,\,   2\al +n \be <\frac2q  +\frac{n}{p} + 1 },\\
D &=  \bket{ (\al, \be) \in (0, 1) \times (0, 1) \,\, : \,\, \be > \frac{n}{(n-1)p}, \,\, 1 +\frac2q < 2\al +\be}.
\end{align*}
Since its verification is straightforward, we skip its details.
\end{rem}

The constructed weak solutions in Theorem \ref{maintheo0503} are not $C^1$ but indeed H\"{o}lder continuous up to the boundary.
Optimal regularity up to boundary is stated in next theorem for the solution of the Stokes system under consideration.

\begin{theo}\label{theo2}
Let   $f$ be given in Assumption \ref{force-f}, and  $ 0 < \al, \, \be < 1$ such that $\ep_0 := 3  -2\al-\be\in  (0, 2) $.  Set $\calR=\bket{(x', x_n ) \in \overline{\mathbb R^n_+} \, : \, |x'| \ge 2,\ x_n \ge 0 }$.
Suppose that $w$ be a solution of of the Stokes system \eqref{maineq}  defined by \eqref{int-0909-1}.
%Let $ A = \{x' \in \Rn \, | \, |x'| >2 \}$.
Then, $w$ is H\"{o}lder continuous in $\calR \times (0,1)$ with the optimal exponent  $\ep_0$, that is,
\begin{equation}\label{0830-1}
 w  \in C^{\ep_0, \frac12  \ep_0}(\calR \times (0,1)),
 \end{equation}
and, in case that $\ep > \ep_0$, we have
\begin{equation}\label{0830-2}
 w \notin  L^\infty (0, 1;  C^{\ep}(\calR) \quad \mbox{and} \quad  w \notin  C^{\frac{\ep}2}(0, 1; L^\infty(\calR)).
\end{equation}
\end{theo}

\begin{rem}
The solution $w$ constructed in (ii) of Theorem \ref{maintheo0503} is contained in  $w \in C^{\ep_0, \frac12  \ep_0}(\calR \times (0,1))$ with $0 < \ep_0 < 1$ and $w \notin C^{\ep, \frac12  \ep}(\calR \times (0,1))$ for all $\ep_0 < \ep$, because $2\al-\be>2$.
On the other hand, in the case that  (ii) of Theorem \ref{maintheo0503-2}, since $\ep_0 := 3 -2\al-\be<2(1-\frac{1}{q})$, $q\in (1, \infty)$, it follows that if $q\le 2$, then $w$  is H\"{o}lder continuous in $\calR \times (0,1)$ with the exponent  $\ep_0\in (0,1)$. On the other hand, in case $q>2$, then $\nabla w$ can be even H\"{o}lder continuous, because $\ep_0$ possibly belongs to $(1,  2(1-\frac{1}{q}))$. We remark that in the interior, the solution $w$ is spatially smooth, although it is just H\"{o}lder continuous in temporal variable (see Proposition \ref{prop0525-3}).
\end{rem}

Lastly, we consider the Navier-Stokes equations in a half space.
\begin{align}\label{NSE-eq}
\left\{\begin{array}{l} \vspace{2mm}
u_t -\De u + {\rm div}\,(u \otimes u)+\na p = f,\quad  {\rm div} \, u =0 \qquad \R_+ \times (0, 1),\\
\vspace{2mm}
\rabs{u}_{x_n =0} =0, \quad \rabs{u}_{t =0} =0.
\end{array}
\right.
\end{align}
Via the method of perturbation, we construct a weak solution of the Navier-Stokes equations whose normal derivatives are unbounded near boundary.
First we specify values of some parameters for such construction.
More pertinently, we choose positive numbers $s$ and $r$  satisfying
\begin{equation}\label{Sept08-10}
\max\bket{\frac{n+2}2, 4} < s < n+2 ,\qquad  \frac{s(n+2)}{n+2-s}< r < \infty.
\end{equation}
Since $2s<\frac{s(n+2)}{n+2-s}$,
it is obvious that $r>2s$. It is also direct  via $n>2$ and \eqref{Sept08-10} that
\begin{equation}\label{Sept08-20}
 2 +\frac{n+2}r< 1  +\frac{n+2}s<2 +\frac{n}2.
\end{equation}
We now fix $\alpha$ and $\beta$ as follows:
\begin{equation}\label{Sept08-30}
\al = 1 - \frac{n+2}{4r} +\frac{\de}{2},\qquad \be = \frac{n+2}{2r} -\ep,
\end{equation}
where $\delta$ and $\ep$ are any number satisfying
$0 <\ep < \de < n \ep < \frac{n+2}{2r}$. It is immediate that $\beta<\frac{1}{2}$, since $r>2s$.
We take $r_0  \in (1,  \infty) $ with $   \frac{3}{r_0} < \de -\ep $.
Let $f$ be the function introduced in Assumption \ref{force-f}  with  $\al$ and $\be$ defined above.
Owing to \eqref{Sept08-10}-\eqref{Sept08-30}, we can see that
\begin{align}\label{Sept08-40}
2 + \frac{3}{r_0} < 2\al + \be, \qquad  2\al + n\be <  2 + \frac{n+2}{r},
\end{align}
and, in addition, reminding, due to \eqref{mixednorm-f}, \eqref{0504-2}, \eqref{maximal-SS} and \eqref{loworder-SS}, it follows for the solution $w$ of  the Stokes system \eqref{maineq} that
\begin{equation}\label{Sept08-50}
w \in L^2(0, 1; \dot H^1 (\R_+)) \cap L^\infty(0,1; L^2 (\R_+)), \quad w \in L^{r} (\R_+\times (0,1)),
\end{equation}
\begin{equation}\label{Sept08-55}
\na w \in L^{s} (\R_+\times (0,1)), \quad  \na w \notin L^{r_0} ((A_i \cup B_i )\times (0,1)\times (0,1)).
\end{equation}

 We are now ready to state the main result for the Navier-Stokes equations.

\begin{theo}\label{thm-NSE}
Let   $f$ be given in Assumption \ref{force-f}, and $\al, \be, s,r$ and  $r_0$ numbers mentioned
in \eqref{Sept08-10}-\eqref{Sept08-40}.
Then there exists a weak solution of the Navier-Stokes equations \eqref{NSE-eq}
so that
\begin{equation}\label{NSE-40}
\na u \in L^{s} (\R_+\times (0,1)),\qquad \na u \notin L^{r_0} ((A_i \cup B_i )\times (0,1)\times (0,1)),
\end{equation}
where $A_i$ and $B_i$ are defined in Notation \ref{boundaryset}.
\end{theo}

\begin{rem}
Here we do not provide definition of weak solutions of the  the Navier-Stokes equations, since it can be
defined similarly as the case of the Stokes system shown in Definition \ref{weaksol-SS}. In such case, the convection term has to be taken into account so that
$\displaystyle \int^1_0\int_{\hR } u\otimes u : \nabla \Phi dxdt$ is to be included in \eqref{weakform-SS}.
\end{rem}

The paper is organized as follows: In Section \ref{section2-0}, we review solution formula of the Stokes system with nonzero force in a half space and prove a lemma that is useful for our analysis.
Section \ref{section2} is devoted to stating a series of propositions that are some parts of main results.
Proofs of propositions are prepared in Section \ref{pfofprops}.
In Section \ref{proofofmaintheo} we present proofs of Theorem \ref{maintheo0503}, Theorem \ref{maintheo0503-2} and Theorem \ref{theo2}. The case of  the Navier-Stokes equations is considered and the proof of Theorem  \ref{thm-NSE} is presented in Section \ref{proof-NSE}.
In Appendix, a figure is drawn to indicate two dimensional cartoon of sets, parts of boundary, where singular solutions of Stokes and Navier-Stokes equations are constructed. In addtion, proofs of technical lemmas such as
Lemma \ref{lemmaappendix} and Lemma \ref{lemma0709-1} are provided.

\section{Preliminaries}
\label{section2-0}
\setcounter{equation}{0}

In this section, we recall the notation of very weak solutions by comparison with weak solutions. We then remind formula of the solution for the Stokes system \eqref{maineq} and some related estimates for the solution. Finally, we provide the proof of  estimates for an integral quantity, which will have beneficial use for our main results.

When given functions $f$ and $g$ are comparable, we use, as a convention of notations,  $f\approx g$, which indicates $c_1 g \le f\le c_2 g$ for some positive constants $c_1$ and $c_2$.

The conception of weak solutions is already introduced and here we account for very weak solutions, which are a bit more generalized weak solutions.

\begin{defin}\label{stokesdefinition}
Let $ f \in L^q_tL^p_x(\mathbb R^n_+ \times (0,1))$ for $1< p,q<\infty$.  A vector field $w\in
L^1_{\rm{loc}}( \hR \times (0, 1))$ is called a very weak solution of the Stokes system \eqref{maineq},
if the following equality is  satisfied:
\begin{align*}
-\int^1_0\int_{\hR }w \cdot \Delta \Phi
dxdt=\int^1_0\int_{\hR}\bke{w \cdot \Phi_t+ f \cdot \Phi} dxdt
\end{align*}
for each $\Phi\in C^2_c(\hR\times [0,1)$ with $\mbox{\rm div } \Phi=0$, and
in addition, for each $\Psi\in C^1_c(\overline{\hR})$
\begin{equation*}
\int_{\hR} w(x,t) \cdot \na \Psi(x) dx =0  \quad  \mbox{ for all}
\quad 0 < t< 1.
\end{equation*}
Furthermore, for each vector field $\varphi\in C^1_c(\overline{\hR})$
\begin{equation*}
\lim_{t\rightarrow 0}\int_{\hR} w(x,t) \cdot \varphi(x) dx=0
\end{equation*}
\end{defin}

For convenience, recalling the formula \eqref{formulas-p}, we decompose the  pressure  $\Pi(x,t) $ as follows:
\begin{equation}\label{formulas-p-1}
\Pi(x,t) =4 \bke{\Pi^{{\mathcal G}}(x,t)+\Pi^{{\mathcal B}}(x,t)},
\end{equation}
where
\begin{equation}\label{CK10-april8}
\Pi^{{\mathcal G}}(x,t)
=   \int_0^t \int_{{\mathbb R}^n_+}
     f_2(y, \tau)
     \int_{{\mathbb R}^{n-1}}
        \Ga(z' -y,  t -\tau) D_{x_2} D_{x_n}
N(x-z')      dz' dyd\tau,
\end{equation}
\begin{equation}\label{badterm-pi}
\Pi^{{\mathcal B}}(x,t)
=  \int_0^t \int_{{\mathbb R}^n_+}
     f_2(y, s)
     \int_{{\mathbb R}^{n-1}} D_{y_n}
        \Ga(z' -y, t -s) D_{x_2}
N (x-z')     dz' dyds.
\end{equation}
For notational conventions, we write the second term of the righthand side in \eqref{formulas-k} as
\begin{equation}\label{denote-L}
L_{ij}(x,y,t):=- 4(1 -\de_{jn}) D_{x_j} \int_0^{x_n} \int_{\Rn}\Ga(x -y^* -z, t) D_{z_i} N(z) dz.
\end{equation}
It was shown in \cite{So1} that  $L_{ij}$  satisfies that for all $ k \in {\mathbb N} \cup \{0\}$, $l= (l', l_n) \in ({\mathbb N} \cup \{0\})^n$,
\begin{align}\label{0801-1}
|D_t^k D_{x_n}^{l_n}D_{x'}^{l'} L_{ij} (x,y,t)| \leq \frac{c e^{-\frac{y_n^2}t}}{ t^k (t + x_n^2)^{\frac{l_n}2} ( |x -y^*|^2 + t )^{\frac{n + |l'| }2}  }, \quad 1 \leq i,j \leq n.
\end{align}

From now on, we denote $Q=\R_+ \times \mathbb (0,1)$, unless any confusion is to be expected.
Next, we recall so called maximal regularity of the Stokes system \eqref{maineq}, which is known to be
as follows: In case that $f\in L^{q_1}_tL_x^{p_1}(Q)$ with $1<p_1, q_1 < \infty, \,\, 1 < p < \infty, \,\, 1 < q \leq \infty$, then
\begin{align}\label{maximal-SS}
& \norm{  w}_{L^{q}_tL_x^{p} (Q)} \le c\norm{f}_{L^{q_1}_tL_x^{p_1} (Q)},\qquad \frac{2}{q}+\frac{n}{p}>\frac{2}{q_1}+\frac{n}{p_1}-2,\\
\label{loworder-SS}
& \norm{D_{x}w}_{L^q_tL_x^p (Q)}\le c\norm{f}_{L^{q_1}_tL_x^{p_1} (Q)},\qquad \frac{2}{q}+\frac{n}{p}> \frac{2}{q_1}+\frac{n}{p_1}-1.
\end{align}

In next lemma we show upper and lower bound estimates of an integral quantity that is related to singular integral of one dimensional heat kernel.

\begin{lemm}\label{lemmaappendix}
Let  $ \be < 1$, $ 0 < \al  < 1$ and $\gamma\in\mathbb R$.  For $x_n >0$ and $t > \frac12 $. We set
\[
\calG(x_n, t):=\int_{\frac12}^t \int_0^2 y_n^{-\be} (\tau-\frac12)^{-\al} (t -\tau)^\ga e^{-\frac{(x_n + y_n)^2}{4(t-\tau)}} dy_n d\tau.
\]
If $\ga -\frac{\be}2 > -\frac32$, then there exist positive constants $c_i$, $i=1,2$ such that
\begin{equation}\label{SS-lemma-10}
c_1 (t-\frac12)^{\frac32 -\frac\be{2}-\al +\ga} e^{-\frac{ x_n^2}{2(t-\frac12)} } \leq \calG(x_n, t) \le c_2 (t-\frac12)^{\frac32 -\frac\be{2}-\al +\ga} e^{-\frac{ x_n^2}{8(t-\frac12)} }.
\end{equation}
If $\ga -\frac{\be}2 \le  -\frac32$, then there exist positive constants $c_i$, $i=3,4$ such that
\begin{equation}\label{SS-lemma-20}
c_3 (t-\frac12)^{ -\al  }x_n^{3 -\be  +2\ga} e^{-\frac{ x_n^2}{2(t-\frac12}) } \leq \calG(x_n, t)  \leq  c_4 (t-\frac12)^{ -\al  }x_n^{3 -\be  +2\ga} e^{-\frac{ x_n^2}{8(t-\frac12)} }.
\end{equation}
\end{lemm}
The proof of Lemma \ref{lemmaappendix} will be presented in Appendix \ref{alter-thm11}.

\section{Stokes equations with external force in a half-space}
\label{section2}
\setcounter{equation}{0}

Let $f$ be the external force defined in {\bf Assumption \ref{force-f}}. For convenience of computations, we decompose $w$ by $w =V + W$, where
\begin{align}\label{V-eqn}
V_i(x,t) &:=  \int_0^t \int_{\R_+}   \big(\Ga(x-y,t-s) - \Ga(x -y^*,t-s) \big) f_i(y,s) dyds,\\
\label{W-eqn}
W_i(x,t)  &:=   \int_0^t \int_{\R_+} L_{i2}(x,y,t-s) f_2(y,s) dyds.
\end{align}
We note that
\begin{equation}\label{0413-1}
D_{x_n} L_{i2}(x,y,t) = D_{x_2} L_{ni}(x,y,t) - 4D_{x_2}  \int_{\Rn} \Ga(x -y^*-z', t) D_{z_i}N(z', 0)d z'.
\end{equation}
Denoting, for convenience,
\begin{align}\label{0706-2}
W_{i}^{{\mathcal G}}(x,t) & =     \int_0^t \int_{\R_+} L_{ni}(x,y,t-\tau) f_2(y,\tau) dyd\tau,\qquad i=1,2,\cdots, n-1,
\end{align}
it follws that
\begin{align}\label{0706-1}
D_{x_n} w_i =D_{x_n} V_i+D_{x_n} W_i =  D_{x_n} V_i +  D_{x_2} W_{i}^{{\mathcal G}}+ {\mathcal B}_{i}^{w},\qquad i=1,2,\cdots, n-1,
\end{align}
 where ${\mathcal B}_{i}^{w}$ is defined as
 \begin{equation}\label{badterm-w}
{\mathcal B}^{w}_i (x,t)
= -4 D_{x_2}\int_0^t \int_{\R_+}   f_2(y,s)  \int_{\Rn} \Ga(x-y^*-z', t-s) D_{z_i}N(z', 0)dz' dyds.
\end{equation}
It turns out that the above term $\calB^w_i$ in \eqref{badterm-w} is the worst term in estimating derivatives for $w$ (see Proposition \ref{theo0521-1}, Proposition \ref{theo0525-1}, Proposition \ref{prop0525-3} and Proposition \ref{prop0525-3-1}).
For the pressure, $\Pi^{{\mathcal B}}$ is more singular than $\Pi^{{\mathcal G}}$ (see Proposition  \ref{prop0503}, Proposition \ref{theo0206} and Proposition \ref{prop0525-3-1}).

On the other hand, since $V$ solves the heat equation in a half space with homogeneous boundary condition, it is worth noting, due to classical regularity theory, that  for $x \in {\mathcal R}$ ($ {\mathcal R}$ is defined in Theorem \ref{theo2}), $ t > \frac12 $ and $i=1,2, \cdots, n$,
\begin{align}\label{0413-3}
\notag | D^l_{x}D_t^m  V_i (x,t) |&  = \abs{  \int_\frac12^t \int_{B_1^+} D^l_{x}D_t^m \Ga(x-y,t-s) f_i (y,s) dyds}\\
\notag &  \leq c   \int_\frac12^t (s -\frac12)^{-\al} (t-s)^{-\frac{n}2 -m -\frac{l}2 } e^{-\frac{c}{t-s}}   ds\,
   \norm{ g }_{L^1 (\R_+)}\\
&  \leq c_k   (t-\frac12)^{k}   \| g \|_{L^1 (\R_+)} \quad \mbox{ for all} \quad k\geq 0.
\end{align}

We start with a key lemma, which is useful in propositions to come after.
Its verification will be given in Appendix \ref{prooflemma0709-1}.
\begin{lemm}\label{lemma0709-1}
Let  $\Ga'$ and $N$ be $n -1$ dimensional Gaussian kernel and  $n$ dimensional Newtonian kernel.
For $x' \neq 0$,  $x_n\ge 0$, $t>0$, $ 1 \leq i \leq n-1$, $k \in ({\mathbb N} \cup \{0\})^{n-1}$ and $  l \geq 0$, If $  |x'| \geq  \max\bket{1, \sqrt{t}}$, it follows that
\begin{align}\label{0515-1}
D^k_{x'} D^l_{x_n}\int_{\Rn}    \Ga'(x'-z',t)    N( z',x_n) dz' =   D^k_{x'} D^l_{x_n}  N( x',x_n)  + J_{kl}(x, t),
\end{align}
such that there exists $c=c(k,l)> 0$, independent of $x$ and $t$, satisfying
\begin{align}\label{Jkl}
\abs{J_{kl}(x, t) } \leq  c
 t^{\frac{1}{2}}.
\end{align}
\end{lemm}

Next, we consider a convolution of second derivatives for Newtonian potential of $N$ and
$g^{{\mathcal T}}$ defined in {\bf Assumption \ref{force-f}}.
More precisely, for $  x \in {\mathcal R}$,  we define
\begin{align}
\label{0830-3}
\phi_{i} (x', x_n) = \int_{|y'| < 1} D_{x_i}D_{ x_2} N(x' -y', x_n) g^{{\mathcal T}}(y') dy'.
\end{align}
We will show that $\phi_{i}$ is strictly positive or negative, depending on the regions under consideration.

Firstly, in case of   $x' \in A_{i1}$ with $i\neq 2$, since $y' \in B'_1$, we note that
\begin{equation}\label{May27-10}
(x_i -y_i) (x_2 -y_2)\geq \frac14 x_i x_2
\geq \frac1{64}( x_i^2 + x_2^2) \geq  \frac1{128}|x'|^2 \geq  \frac1{512}|x' -y'|^2.
\end{equation}
Conversely, if $x' \in A_{i2}$, then we can see that
\begin{equation}\label{May27-20}
(x_i -y_i) (x_2 -y_2) \leq \frac14 x_i x_2
\leq  - \frac1{64}( x_i^2 + x_2^2) \leq  -\frac1{128}|x'|^2 \leq  -\frac1{512}|x' -y'|^2.
\end{equation}
On the other hand, in case that $x' \in B_{i1}$ and $y' \in B'_1$, it follows that
\begin{equation}\label{May27-30}
|x'-y'|^2-n(x_2-y_2)^2  \geq \frac1{64} |x' -y'|^2.
\end{equation}
Indeed, recalling \eqref{Bi1-set}, it is straightforward that
\[
|x'-y'|^2-n(x_2-y_2)^2  \ge\frac{1}{4}|x'|^2-2n|x_2|^2\ge \frac{1}{16}|x'|^2
\ge \frac1{64} |x' -y'|^2.
\]
Similarly, we can see that if $ x' \in B_{i2}$, then
\[
|x'-y'|^2-n(x_2-y_2)^2 \leq  - \frac1{64} |x' -y'|^2.
\]
Hence, for  $x' \in A_{i1}$ with $ i \neq 2$,  we observe that
\begin{align}\label{est-Ai1}
\nonumber \phi_i (x',x_n)& = -c \int_{\Rn} g^{{\mathcal T}}(y')\frac{(x_i -y_i) (x_2 -y_2) }{|x-y'|^{n+2}} dy'\\
& \leq -c  \int_{\Rn} g^{{\mathcal T}}(y')\frac{1}{|x -y'|^n} dy', %\quad x' \in A_{i1},
\end{align}
where we used \eqref{May27-10}. Analogously, in case that $x' \in A_{i2}$ with $ i \neq 2$, it follows that
\begin{align}\label{est-Ai2}
\nonumber \phi_i (x', x_n)& = -c\int_{\Rn} g^{{\mathcal T}}(y') \frac{(x_i -y_i) (x_2 -y_2) }{|x-y'|^{n+2}}dy_n\\
&  \geq c\int_{\Rn} g^{{\mathcal T}}(y')\frac{1}{|x -y'|^n} dy',
\end{align}
where \eqref{May27-20} is used.
Meanwhile, for $x' \in B_{i1}$  with $ i \neq 2$ we also note via \eqref{May27-30} that
\begin{align}\label{est-Bi1}
\nonumber \phi_2 (x', x_n)& = c \int_{\Rn} g^{{\mathcal T}}(y') \frac{|x-y'|^{-2}-n(x_2-y_2)^2}{|x-y'|^{n+2}} dy'\\
& \geq c\int_{\Rn} g^{{\mathcal T}}(y') \frac{1}{|x -y'|^n}dy'. %\quad x' \in B_{i1},
\end{align}
Likewisely, for  $x' \in B_{i2}$ with $ i \neq 2$ we observe that
\begin{align}\label{est-Bi2}
\nonumber \phi_2(x', x_n) &  = -c \int_{\Rn} g^{{\mathcal T}}(y') \frac{|x-y'|^{-2}-n(x_2-y_2)^2 }{|x-y'|^{n+2}}dy'\\
&
\leq - c \int_{\Rn} g^{{\mathcal T}}(y')\frac{1}{|x-y'|^n}dy'.
\end{align}

Next proposition shows pointwise estimates of the worst term ${\mathcal B}^{w}_i$ defined in \eqref{badterm-w}. In particular, with aid of \eqref{est-Ai1}-\eqref{est-Bi2},
lower bounds or upper bounds are provided on mutually disjoint sets near boundary. All proofs of propositions in this section will be provided in Section \ref{pfofprops}.

\begin{prop}\label{theo0521-1}
Let $1 \leq i \leq n-1$. Suppose  that  ${\mathcal B}^{w}$ is defined in \eqref{badterm-w} and $\phi_i$ is defined in \eqref{0830-3}. Then, for $x \in {\mathcal R} $ and $t>\frac{1}{2}$
\begin{align}\label{0522-1}
\abs{{\mathcal B}^{w}_i  (x,t)}\,\,
\ge \left\{\begin{array}{l}
c (t -\frac12)^{1 -\frac{\be}2 -\al} e^{-\frac{x_n^2}{2(t-\frac12)}} \phi_i(x',0)\chi_{(A_{i1}\cup B_{i1})}+O\bke{(t-\frac12)^{\frac32 -\frac{\be}2 -\al }},\\
 c (t -\frac12)^{1 -\frac{\be}2 -\al} e^{-\frac{x_n^2}{8(t-\frac12)}} \phi_i(x',0)\chi_{(A_{i2}\cup B_{i2})}+O\bke{(t-\frac12)^{\frac32 -\frac{\be}2 -\al }}.
\end{array}
\right.
\end{align}
More precisely,  for $i \neq 2$ and $t > \frac12$,
\begin{align}\label{0522-2}
{\mathcal B}^{w}_i  (x,t) & \leq   -c (t-\frac12)^{1-\frac{\be}2 -\al } e^{-\frac{x_n^2}{2(t-\frac12)}}\phi_i(x', 0)+ c(t-\frac12)^{\frac32 -\frac{\be}2 -\al } , \quad  x' \in A_{i1},\\
\label{0522-3}{\mathcal B}^{w}_i  (x,t) & \geq   c (t-\frac12)^{1 -\frac{\be}2 -\al } e^{-\frac{x_n^2}{8(t-\frac12)}}\phi_i(x', 0)- c(t-\frac12)^{\frac32 -\frac{\be}2 -\al } , \quad   x' \in A_{i2},\\
\label{0522-4}
{\mathcal B}^{w}_2 (x,t) & \leq   -c (t-\frac12)^{1 -\frac{\be}2 -\al } e^{-\frac{x_n^2}{2(t-\frac12)}} \phi_2(x', 0)+ c(t-\frac12)^{\frac32 -\frac{\be}2 -\al } , \quad   x' \in B_{i1},\\
\label{0522-5}{\mathcal B}^{w}_2 (x,t) & \geq   c (t-\frac12)^{1 -\frac{\be}2 -\al } e^{-\frac{x_n^2}{8(t-\frac12)}}\phi_2(x', 0) - c(t-\frac12)^{\frac32 -\frac{\be}2 -\al } , \quad  x' \in B_{i2},
\end{align}
where $A_{i1}, A_{i2}, B_{i1}$ and $B_{i2}$ are sets introduced  in {\bf Notation \ref{boundaryset}}.
\end{prop}

We can also obtain estimates of higher order derivatives for the term ${\mathcal B}^{w}_i$.

\begin{prop}\label{theo0525-1}
Suppose  that  ${\mathcal B}^{w}$ is defined in
\eqref{badterm-w}.
Let $1 \leq i \leq n-1$.  Then, for $x \in {\mathcal R} $ and $t>\frac{1}{2}$
\begin{align}\label{0706-7}
| D_{x'}^k D_{x_n}^l   {\mathcal B}^{w}_i
(x,t) | & \leq    c \left\{\begin{array}{l}
(t-\frac{1}{2})^{1-\frac{1+\be+2\al}2 }e^{-\frac{x_n^2}{8(t-\frac{1}{2})}}, \quad l  = 1\\
(t-\frac12)^{ -\al  }x_n^{2 -\be  -l} e^{-\frac{ x_n^2}{8(t-\frac12)} }, \quad l \geq 2,
\end{array}
\right.
\end{align}
where $c_2$ and $c_4$ are constants in Lemma \ref{lemmaappendix}.
For a given positive integer $l$, there exists $c_l>0$ such that  if  $\sqrt{t -\frac12} \leq c_l x_n$, then
\begin{align}\label{0706-8}
| D_{x'}^k D_{x_n}^l   {\mathcal B}^{w}_i
(x,t) | & \geq   c
(t-\frac{1}{2})^{1-\frac{l+\be+2\al}2 }e^{-\frac{x_n^2}{2(t-\frac{1}{2})}}.
\end{align}
\end{prop}

We haven't mentioned any integrability of the pressure corresponding to the weak solution so far. In next proposition, we show that the pressure belongs to some Lebesgue spaces, in particular, globally in space.

\begin{prop}\label{prop0503}
Suppose  that $\Pi$  is given in \eqref{int-expression-p} and
$f\in L^{q_1} (0, 1; L^{p_1} (\R_+))$ with $1<p_1, q_1<\infty$.
If $p$ and $q$ satisfy
$p>\frac{n}{n-1}p_1 $ and $\frac2{q_1}+\frac{n}{p_1}-1 \leq  \frac2q  + \frac{n}{p}$,
then
\begin{equation}\label{pi-bound-2}
\| \Pi\|_{L^{q}(0, 1; L^p (\R_+))} < c \| f_2\|_{L^{q_1} (0, 1; L^{p_1} (\R_+))}.
\end{equation}
\end{prop}

Next proposition shows pointwise estimates of spatial derivatives for $\Pi^{{\mathcal B}}$
 near boundary.
 For simplicity, ${\rm sgn(a)}$ indicates $1$, if $a>0$ and otherwise, ${\rm sgn(a)}=-1$, and we also denote
\begin{equation}\label{june16-100}
\psi (x) := \int_{\Rn}   D_{x_2}   N( x'-y',x_n) g^{{\mathcal T}} (y')   dy'.
\end{equation}
\begin{prop}\label{theo0206}
Suppose $\Pi^{{\mathcal B}}$ is defined in
\eqref{badterm-pi}.
Then, for $|x_2| \geq 2$ and $ \frac12 < t \leq 1$,
\begin{align}\label{0830-4}
\abs{\abs{D_{x'}^kD_{x_n}^l\Pi^{{\mathcal B}}(x,t)}- c(t-\frac12)^{\frac12 -\frac{\be}2 -\al} \abs{D_{x'}^kD_{x_n}^l\psi (x)}}\le c(t-\frac12)^{1 -\frac{\be}2 -\al},\qquad k \geq 0.
\end{align}
In particular, in case that $k=l=0$, we have
\begin{align}\label{june17-10}
\abs{\Pi^{{\mathcal B}}(x,t) +c \,{\rm sgn}(x_2)\,(t-\frac12)^{\frac12 -\frac{\be}2 -\al} \psi(x)}\le c(t-\frac12)^{1 -\frac{\be}2 -\al},
\end{align}
where $\chi_A$ is the characteristic function supported on $A$ and $\psi$ is defined in \eqref{june16-100}.
\end{prop}

\begin{rem}
Let $ 1 < p_1  < \frac{n-1}np < \infty$ and $1 < q_1<  q < \infty$. We denote
\begin{align*}
C &= \{ (\al, \be) \in (0, 1) \times (0, 1) \,| \,  \frac1{q_1} -\frac1{q} +\frac{n}{2p_1} -\frac{n}{2p} \leq \frac12, \,\, \al < \frac1{q_1}, \,\, \be < \frac1{p_1} \},\\
D&  =  \{ (\al, \be) \in (0, 1) \times (0, 1) \, | \, \frac12 + \frac1{q}     < \al +\frac{\be}2 \}.
\end{align*}
Then, $C \cap D =\emptyset$, which indicates that \eqref{june17-10} is not conflict to \eqref{pi-bound-2}.
\end{rem}

In next proposition, we state some estimates of
velocity fields and the pressure excluding  $\calB^w_i$ and $\Pi^{{\mathcal B}}$, respectively. Those estimates turn out to be less singular than the worst terms $\calB^w_i$ and $\Pi^{{\mathcal B}}$, as mentioned earlier.

\begin{prop}\label{prop0525-3}
Let $1 \leq i \leq n-1$. Suppose  that $w$ and $\Pi$  are given in \eqref{int-0909-1} and \eqref{int-expression-p}, and ${\mathcal B}^{w}$ and $\Pi^{{\mathcal B}}$ are defined in
\eqref{badterm-w} and \eqref{badterm-pi}, respectively. Then, for any  $ k, \, l \geq 0$ and $x \in {\mathcal R} $, $t>\frac{1}{2}$,  we have
\begin{align}\label{0525-2}
\abs{ D^{k}_{x'} D_{x_n}^{l}\bke{D_{x_n} w_i (x,t) - {\mathcal B}^{w}_i  (x,t)}} & \le  c\left\{\begin{array}{l}
(t-\frac{1}{2})^{1-\frac{l-2+\be+2\al}2 }e^{-\frac{x_n^2}{8(t-\frac{1}{2})}}, \quad  1 \leq l \leq 3,\\
(t-\frac12)^{ -\al  }x_n^{2 -\be  -l} e^{-\frac{ x_n^2}{8(t-\frac12)} }, \quad l \geq 4,
\end{array}
\right.
\end{align}
\begin{equation}\label{0616-10}
\abs{ D^{k}_{x'} D_{x_n}^{l+1}w_n(x,t)}\le c\left\{\begin{array}{l}
(t-\frac{1}{2})^{1-\frac{l-1+\be+2\al}2 }e^{-\frac{x_n^2}{8(t-\frac{1}{2})}}, \quad l \leq 1,\\
(t-\frac12)^{ -\al  }x_n^{3 -\be  -l} e^{-\frac{ x_n^2}{8(t-\frac12)} }, \quad l \geq 2,
\end{array}
\right.
\end{equation}
and
\begin{align}\label{0525-1}
\abs{D^k_{x}\Big(\Pi (x,t) - \Pi^{{\mathcal B}}(x,t) \Big) } \leq c t^{1 -\frac{\be}2 -\al}, \qquad k\ge 0.
\end{align}
 \end{prop}

Next proposition shows  some estimates involving temporal derivative of velocity.

\begin{prop}\label{prop0525-3-1}
Suppose  that $w$ and $\Pi$  are given in \eqref{int-0909-1} and \eqref{int-expression-p}, and ${\mathcal B}^{w}$ and $\Pi^{{\mathcal B}}$ are defined in
\eqref{badterm-w} and \eqref{badterm-pi}, respectively. Then, for $1 \leq i \leq n-1$ and  for any $ k, l\ge 0$ and $x \in {\mathcal R}, t>\frac{1}{2}$,
\begin{equation*}
\abs{D_{x'}^kD^l_{x_n} \bke{  D_t w_i (x,t) -  D_{x_n} {\mathcal B}^{w}_i (x,t)+D_{x_i}\Pi^{{\mathcal B}}(x,t) }  }
\end{equation*}
\begin{equation}\label{May27-90}
\leq c \max\bket{(t-\frac12)^{ -\al  }x_n^{2 -\be  -l} e^{-\frac{ x_n^2}{8(t-\frac12)} }, (t-\frac{1}{2})^{1 -\frac{\be}2 -\al}}.
\end{equation}
The normal component $w_n$ satisfies
\[
\abs{ D^{k}_{x'} D_{x_n}^{l+1}\bke{ D_t w_n (x,t) -D_{x_n}^{2} w_n(x,t)+ D_{x_n}\Pi^{{\mathcal B}}(x,t) }}
\]
\begin{equation}\label{0617-100}
\le c\max\bket{(t-\frac12)^{ -\al  }x_n^{2 -\be  -l} e^{-\frac{ x_n^2}{8(t-\frac12)} },  (t-\frac{1}{2})^{1 -\frac{\be}2 -\al}},
\end{equation}
\end{prop}

\begin{rem}
The meaning of \eqref{May27-90} is the control of $D_t w_i -\Delta w_i-D_{x_i} \Pi$, $i\neq n$ near boundary is majorized by
$D_t w_i  -  D_{x_n} {\mathcal B}^{w}_i +D_{x_i}\Pi^{{\mathcal B}}$.
Similarly, the main profile of $D_t w_n -\Delta w_n+D_{x_n}  \Pi$ is $D_t w_n  -D_{x_n}^{2} w_n+ D_{x_n}\Pi^{{\mathcal B}}$ for  the normal component of the equation.
\end{rem}

\section{Proofs of Propositions}\label{pfofprops}
\setcounter{equation}{0}

\subsection{Proof of Proposition \ref{theo0521-1}}

\label{prprop0521-1}

Note that since $|x'| \geq 2$, it follows that $|x' -y'| \geq 1$ for $y' \in {\rm supp} \, g^{{\mathcal T}} \subset B_1^{'}$. Using   \eqref{0515-1} in Lemma \ref{lemma0709-1}, we decompose ${\mathcal B}^{w}_i $ as follows:
\begin{align}\label{0710-3}
 \notag  {\mathcal B}^{w}_i (x,t) &=   -c_n \int_{0}^t \int_{\R_+} f_2(y, \tau)   \frac{1}{(t -\tau)^{\frac{1}2}}e^{-\frac{(x_n + y_n)^2}{4(t-\tau)}}\\
\notag  & \quad
\times D_{x_2}  \int_{\Rn} \Ga'(x'-y' -z',t-\tau) D_{z_i}N(z', 0)  dz' dy d\tau\\
& : = I_{i1}(x,t) + I_{i2}(x,t),
\end{align}
where
\begin{align*}
I_{i1} (x,t)& = -c_n \int_{0}^t  \int_{\R_+}  f_2(y,\tau)   \frac{1}{(t -\tau)^{\frac{1}2}}e^{-\frac{(x_n + y_n)^2}{4(t-\tau)}} D_{x_2} D_{x_i} N(x' -y',0)
 dy d\tau,\\
I_{i2}(x,t)& = -c_n \int_{0}^t \int_{\R_+}    f_2(y,\tau)   \frac{1}{(t -\tau)^{\frac{1}2}}e^{-\frac{( x_n +y_n)^2}{4(t-\tau)}} J_{10}(x' -y', t-\tau)
 dy d\tau.
\end{align*}
By  \eqref{Jkl}, we have
\begin{align*}
\int_{ \Rn}g^{{\mathcal T}}(y') J_{10}(x' -y', t-\tau) dy' & \leq   \int_{|y'| \leq 1 }  g^{{\mathcal T}}(y') (t -\tau)^\frac12  dy'\\
& \leq \|g^{{\mathcal T}}\|_{L^\infty} (t -\tau)^{\frac12}.
\end{align*}
Note that $|g^{\mathcal N}(y_n)| \leq c y^{-\al}_n$ for $0< y_n < 2$. From Lemma \ref{lemmaappendix},for $x \in {\mathcal R} $ and $t>\frac{1}{2}$,  we obtain
\begin{align}
\notag \|I_{i2} (t) \|_{L^\infty(|x'| >2 )}  & \leq c  \| g^{{\mathcal T}}\|_{L^\infty(\Rn)}\int_{\frac12}^t \int_0^2   (\tau -\frac12)^{-\al} y_n^{-\be}     e^{-\frac{ (x_n + y_n)^2}{8(t-\tau)}} dy_n d\tau\\
\label{0413-5}& \leq  c  \| g^{{\mathcal T}}\|_{L^\infty(\Rn)} (t-\frac12)^{\frac32 -\frac{\be}2 -\al} e^{-\frac{ x_n^2}{8(t-\frac12)}}.
\end{align}
By Lemma \ref{lemmaappendix}, for $|x'| \geq 2  $, we have
\begin{align}\label{0413-6}
\notag I_{i1}(x,t) &  =- c \int_{\frac12}^t \int_0^2   (\tau -\frac12)^{-\al} y_n^{-\be}   (t -\tau)^{-\frac12 }   e^{-\frac{ (x_n + y_n)^2}{4(t-\tau)}} dy_n d\tau \phi_i(x',0)\\
& \left\{\begin{array}{l}
\le  -c (t -\frac12)^{1 -\frac{\be}2 -\al} e^{-\frac{x_n^2}{2(t-\frac12)}} \phi_i(x',0),\\
\ge -c (t -\frac12)^{1 -\frac{\be}2 -\al} e^{-\frac{x_n^2}{8(t-\frac12)}} \phi_i(x',0).
\end{array}
\right.
\end{align}
Summing all estimates, we obtain \eqref{0522-1} by taking the absolute value.
Using from \eqref{est-Ai1} to \eqref{est-Bi2}, we get from  \eqref{0522-2}-\eqref{0522-5}.
We complete the proof of Proposition  \ref{theo0521-1}.

\subsection{Proof of Proposition \ref{theo0525-1}}
\label{prooftheo0206}
Let $\Gamma_1(x_n,t)$ be the heat kernel in one dimension. For simplicity, we denote $\eta=\frac{x_n}{2\sqrt{t}}$. We then note that
\begin{equation}\label{June19-10}
D_{x_n}^l \Gamma_1(x_n, t)=\frac{1}{2\sqrt{\pi}}t^{-\frac{l+1}{2}}P_l(\eta)e^{-\eta^2},
\qquad P_l(\eta)=\sum_{i=0}^{\bkt{\frac{l}{2}}}c_{li}\eta^{2(i +\frac{l}{2}-\bkt{\frac{l}{2}})},
\end{equation}
where $c_{li}$ is a nonzero constant.
For example,
\[
P_0(\eta)=1, \quad P_1(\eta)=-2\eta,\quad P_2(\eta)=-2+4\eta^2,\quad P_3(\eta)=12\eta-8\eta^3.
\]
In fact, $P_l$ satisfies the following recursive formula:
\[
P_{l}(\eta)=P'_{l-1}(\eta)-2\eta P_{l-1}(\eta), \quad P_0 (\eta)=1,\qquad l\ge 1.
\]

Again, we only compute normal derivatives, since the tangential derivatives are rather easy. With aid of Lemma \ref{lemma0709-1}, for $x \in {\mathcal R} $ and $t>\frac{1}{2}$, we have

\begin{align}\label{0522-6}
\nonumber D^l_{x_n} {\mathcal B}^{w}_i (x,t)
 = & -4  \int_\frac12^t \int_0^2    g^{\mathcal N}(y_n) h(s)    D^l_{x_n} \Ga_1(x_n +  y_n, t-\tau) dy_n d\tau \phi_i(x',0)\\
\nonumber&   -4 \int_\frac12^t \int_0^2    g^{\mathcal N}(y_n) h(s)   D^l_{x_n} \Ga_1(x_n +  y_n, t-\tau) J_{0l}(x', 0, t-\tau) dy_nd\tau \\
 = & I_1 + I_2.
\end{align}
Using \eqref{June19-10}, by Lemma \ref{lemmaappendix},   we obtain
\begin{align}\label{June20-10}
 \nonumber |I_2| & \leq c \int_\frac12^t  (\tau-\frac{1}{2})^{-\al} (t -\tau)^{-\frac{l}{2}}   \int_0^2   y_n^{-\be}   e^{-\frac{(x_n +y_n)^2}{ 8(t-\tau)}} dy_nds\\
 & \leq c \left\{\begin{array}{ll} (t-\frac{1}{2})^{\frac{3-l-\be-2\al}2 }e^{-\frac{x_n^2}{8(t-\frac{1}{2})}}, \quad l \leq 2\\
 (t-\frac12)^{ -\al  }x_n^{3 -\be  -l} e^{-\frac{ x_n^2}{8(t-\frac12)} }, \quad l \geq 3,
 \end{array}
 \right.
\end{align}
where we used that $x_n^2 + y_n^2 \leq  (x_n + y_n)^2 \leq 2 (x_n^2 + y_n^2)$ for $x_n, \,\, y_n \geq 0$.
On the other hand, to estimate, by Lemma \ref{lemmaappendix}, we compute
\begin{align}\label{June20-20}
\notag |I_1| & =c\abs{ \int_0^t \int_0^2 y_n^{-\be}  (\tau-\frac{1}{2})^{-\al}     D^l_{x_n} \Ga_1(x_n +  y_n, t-\tau) dy_n d\tau }  |\phi_i(x',0)| \\
\notag & \leq  c\int_0^t   (\tau-\frac{1}{2})^{-\al}   (t -\tau)^{-\frac{1+l}2}     \int_0^1 y_n^{-\be}    e^{-\frac{ (x_n +y_n)^2}{2(t-\tau)}} dy_n d\tau  |\phi_i(x',0)|\\
& \leq  c\left\{\begin{array}{l}
(t-\frac{1}{2})^{1-\frac{l+\be+2\al}2 }e^{-\frac{x_n^2}{8(t-\frac{1}{2})}} |\phi_i(x',0)|, \quad l \leq 1\\
(t-\frac12)^{ -\al  }x_n^{2 -\be  -l} e^{-\frac{ x_n^2}{8(t-\frac12)} }|\phi_i(x',0)|, \quad l \geq 2.
\end{array}
\right.
\end{align}
Summing up the above estimates \eqref{June20-10} and \eqref{June20-20}, we obtain \eqref{0706-7}.

For the lower bounded, note that there are $c_{l}> 0$ such that if  $\eta > c_{l}$ then
\begin{align}\label{1014-10}
|P_{l} (\eta)| \geq \eta^l.
\end{align}
Let $ t-\frac12 < c_l^{-2} x_n^2$. We note that if    $\frac12 < \tau$, then $c_l \sqrt{t -\tau} -x_n < 0$.
Hence, due to \eqref{1014-10}, we have
\begin{align}\label{0706-6}
\notag |I_1| & \geq
  c\int_\frac12^t   (\tau-\frac{1}{2})^{-\al}   (t -\tau)^{-\frac{1+l}2}     \int_0^1 y_n^{-\be} (\frac{ x_n +y_n}{t-\tau})^l   e^{-\frac{ (x_n +y_n)^2}{4(t-\tau)}} dy_n d\tau  |\phi_i(x',0)|\\
\notag   & \geq
  c \int_\frac12^t   (\tau-\frac{1}{2})^{-\al}   (t -\tau)^{-\frac{1+l}2}     e^{-\frac{ x_n^2}{2(t-\tau)}}  \int_0^1 y_n^{-\be} (\frac{ y_n}{t-\tau})^l   e^{-\frac{ y_n^2}{2(t-\tau)}} dy_n d\tau  |\phi_i(x',0)|\\
\notag   & =
  c  \int_\frac12^t   (\tau-\frac{1}{2})^{-\al}   (t -\tau)^{-\frac{1+l}2 +\frac12 -\frac{\be}2 }     e^{-\frac{ x_n^2}{2(t-\tau)}}  \int_0^{\frac1{\sqrt{t -\tau}}} y_n^{l-\be}     e^{-  y_n^2 } dy_n d\tau  |\phi_i(x',0)|\\
& \geq  c
(t-\frac{1}{2})^{1-\frac{l+\be+2\al}2 }e^{-\frac{x_n^2}{2(t-\frac{1}{2})}} |\phi_i(x', 0)|.
\end{align}
By  \eqref{June20-10} and \eqref{0706-6},  we obtain \eqref{0706-8}. We complete the proof of Proposition \ref{theo0525-1}.

\subsection{Proof of Proposition \ref{prop0503}}
\label{proofprop0503}

Note that $ \Pi = \Pi^{{\mathcal G}} + \Pi^{{\mathcal B}}$, where
\begin{align*}
\Pi^{{\mathcal G}}(x,t) = P_{n}   \bke{D_{x_2} \rabs{\Ga * f_2}_{x_n =0}}(x',t), \qquad
\Pi^{{\mathcal B}}(x,t) = P_{n} R'_2 \bke{D_{x_n} \rabs{\Ga * f_2}_{x_n =0}}(x',t),
\end{align*}
where $P_n$ is the Poisson Kernel of the Laplace equation in $\mathbb R^n_+$ and $R'_2$ is the Riesz transform in $\mathbb R^{n-1}$ for $x_2$ variable.
Let $ 0<\ep < 1$ be sufficiently  small positive number and $\frac1{r_1} = \frac{n}{p(n-1)} +\ep$  satisfying  $p_1 < r_1 < \frac{n-1}n p$. Take $\ep > 0$ small such that  $  \te: = \frac{n}{r_1} -\frac{n}p = \frac{n}{p (n-1)} +\ep n  = \frac1{r_1} +\ep (n-1)$ satisfies $ 0 <  \te < 1$. Choose $1< r_2, \,\, q_2 < \infty$ satisfying $\frac1{q_1} = \frac{1-\te}{q_2} +\frac{\te}{q}$ and $\frac1{p_1} = \frac{1-\te}{r_2} +\frac{\te}{r_1}$. If $\ep > 0$ is sufficiently small, then $1 < r_2 < p_1$ and $1 < q_2 < q_1 $. Note that $\frac{n-1}{r_1} -\frac{n}p >0$. From well-known harmonic function estimate, Besov inequality and trace theorem,   we have
\begin{align*}
\| \Pi (t) \|_{L^p (\R_+)}
& \le c \big(  \|D_{x_n} \Ga * f_2(t)|_{x_n =0}\|_{\dot B^{-\frac1p}_p (\Rn)} + \|D_{x_2} \Ga * f_2(t)|_{x_n =0}\|_{\dot B^{-\frac1p}_p (\Rn)} \big) \\
& \le c   \|D_{x} \Ga * f_2(t)|_{x_n =0}\|_{\dot B^{\frac{n-1}{r_1} -\frac{n}p}_{r_1} (\Rn)} \\
& \le c \|D_{x} \Ga * f_2(t) \|_{\dot H^{\frac{n}{r_1} -\frac{n}p}_{r_1} (\R_+)}.
\end{align*}
Hence, we have
\begin{align*}
\| \Pi  \|_{L^q(0, \infty; L^p (\R_+))}
& \le c \|D_{x} \Ga * f_2 \|_{L^q (0, \infty;\dot H^{\frac{n}{r_1} -\frac{n}p}_{r_1} (\R_+))}.
\end{align*}

Note that $\frac1q +1 = \frac1{q_2} +\frac12 +\frac{n}{2r_2} -\frac{n}{2 r_1}$. Then, we have
\begin{align*}
\|D_{x} \Ga * f_2 \|_{L^q (0, \infty;  L^{r_1} (\R_+))}
& \leq c \| f_2 \|_{L^{q_2} (0, \infty;  L^{r_2} (\R_+))},\\
\|D_{x} \Ga * f_2 \|_{L^q (0, \infty; \dot H^1_{r_1} (\R_+))}
& \leq c \| f_2 \|_{L^{q} (0, \infty;  L^{r_1} (\R_+))}.
\end{align*}
Using the complex interpolation property, we have  $[L^q (0, \infty; L^{r_1} (\R_+)),  L^q (0, \infty;\dot H^1_{r_1} (\R_+))]_\te = L^q (0, \infty;\dot H^\te_{r_1} (\R_+))$ and
  $[L^{q_2} (0, \infty; L^{r_2} (\R_+), L^q (0, \infty; \dot L^{r_1} (\R_+)]_\te = L^{q_1} (0, \infty; L^{p_1} (\R_+)) $ (see \cite{BL}).  Hence, we have
\begin{align*}
\| \Pi  \|_{L^q(0, \infty; L^p (\R_+))}
& \le c \|D_{x } \Ga * f_2 \|_{L^q (0, \infty;\dot H^{\te}_{r_1} (\R_+))} \leq c \| f_2 \|_{L^{q_1} (0, \infty; L^{p_1} (\R_+))}.
\end{align*}
Hence, we complete the proof of Proposition \ref{prop0503}.

\subsection{Proof of Proposition \ref{theo0206}}
\label{prooftheo0206}

With the aid of \eqref{0515-1}, we decompose $\Pi^{{\mathcal B}}$ as follows:
\begin{align*}
D_{x'}^kD_{x_n}^l \Pi^{{\mathcal B}}(x,t) = &- 4 c_n \int_0^t \int_{{\mathbb R}^n_+}
     f_2(y, \tau) (t-\tau)^{-\frac32}y_n    e^{-\frac{y_n^2}{t-\tau}}
      D_{x'}^kD_{x_n}^l D_{x_2}    N( x'-y',x_n)  dyd\tau \\
  &- 4 c_n \int_0^t \int_{{\mathbb R}^n_+}
     f_2(y, \tau) (t-\tau)^{-\frac32}y_n    e^{-\frac{y_n^2}{t-\tau}} J_{k+1,l}(x'-y', t-\tau)  dyd\tau\\
 :=& I+J.
\end{align*}
Since for $|x_2| \geq 2$, $ t > \frac12$ and $y' \in B'_1$, by \eqref{Jkl} and \eqref{SS-lemma-10}, we have
\begin{align*}
|J|
  &  \leq c \int_{\frac{1}{2}}^t  (\tau-\frac{1}{2})^{-\al}     (t-\tau)^{-\frac{\be}2}     d\tau
   \leq c (t-\frac12 )^{1 -\frac{\be}2 -\al}.
\end{align*}
Recalling \eqref{june16-100}, we note  for $|x'| >2$, by \eqref{SS-lemma-10} that
\begin{align*}
I& = -c \int_0^t  h(\tau) \int_0^2 g^{\mathcal N}(y_n) y_n     (t-\tau)^{-\frac32}    e^{-\frac{y_n^2}{t-\tau}}
    dy_n d\tau   D_{x'}^kD_{x_n}^l \psi(x', x_n)  \\
  & \approx -  (t-\frac12 )^{\frac12 -\frac{\be}2 -\al}   D_{x'}^kD_{x_n}^l \psi(x', x_n),
\end{align*}
which implies \eqref{0830-4} and  \eqref{june17-10}. Therefore,
we complete the proof of  Proposition \ref{theo0206}.

\subsection{Proof of Proposition \ref{prop0525-3}}
\label{prprop0525-3}

We begin with spatial estimates of $D_{x_n} w_i (x,t) - {\mathcal B}^{w}_i  (x,t)$.
Due to \eqref{0706-1} and \eqref{0413-3} , it is enough to estimate spatial derivatives of $D_{x_2} W_{i}^{{\mathcal G}}$.

$\bullet$ (case $l=0$)\,\, It suffices to estimate the $k$-th order tangential derivatives of $D_{x_2} W_{i}^{{\mathcal G}}$. Recalling the {\bf Assumption \ref{force-f}} and using the
estimates \eqref{0801-1} and \eqref{SS-lemma-10}, it follows for $x \in {\mathcal R} $ and $t>\frac{1}{2}$ that
\begin{equation}\label{L0-10}
\abs{ D_{x'}^k D_{x_2} W_{i}^{{\mathcal G}}(x,t)}\le c \int_\frac12^t \int_{0}^2 y_n^{-\beta}(\tau-\frac{1}{2})^{-\alpha} e^{-\frac{y^2_n}{t-\tau }}dy_nd\tau\le c(t-\frac{1}{2})^{\frac32  -\frac{\be}{2} -\al}.
\end{equation}

Thus, combining \eqref{0706-1},  \eqref{0413-3} and \eqref{L0-10}, we obtain
\[
\abs{ \bke{D_{x_n} w_i (x,t) - {\mathcal B}^{w}_i  (x,t)}}\le c(t-\frac{1}{2})^{\frac32  -\frac{\be}{2} -\al}.
\]

$\bullet$ (case $l=1$)\,\,
We note first that the normal derivatives bring changes in estimates and
tangential derivative doesn't make any difference, and therefore, we compute only normal derivatives.
On the other hand, we observe that
\begin{align}\label{june16-10}
\notag D_{x_n} L_{ni}(x,t) & =
4D_{x_n} D_{x_i} \int_0^{x_n} \int_{\Rn} \Ga(x -y^* -z, t) D_{z_n} N(z) dz\\
& =  -\sum_{ m=1}^{ n-1} D_{x_i} L_{mm}(x,y,t) + 2D_{x_i}\Ga(x -y^*, t),
\end{align}
which implies again via  \eqref{0706-2} that
\begin{align*}
D_{x_n} D_{x_2} W_i^{{\mathcal G}} (x,t) =  -\sum_{ 1 \leq m \leq n-1}D_{x_2}  D_{x_i} W_{mm} (x,t) + 2 D_{x_2} D_{x_i} \Ga^* * f_2(x,t),
\end{align*}
where
\begin{align*}
\Ga^* * f_2(x,t)&=\int_0^t \int_{\R_+} \Ga(x -y^*, t-\tau) f_2 (y,\tau) dyd\tau, \\
W_{m_1m_2} (x,t) &= \int_0^t \int_{\R_+} L_{m_1m_2} (x,y,t-\tau) f_2 (y,\tau) dyd\tau, \quad 1 \leq m_1, m_2 \leq n.
\end{align*}
To sum up, we obtain
\begin{align}\label{0706-3}
D_{x_n} \big( D_{x_n} w_i -   {\mathcal B}_{i}^{w} \big)  = D^2_{x_n} V_i -\sum_{ 1 \leq m \leq n-1}D_{x_2}  D_{x_i} W_{mm} (x,t) + 2 D_{x_2} D_{x_i} \Ga^* * f_2(x,t).
\end{align}
Repeating similar computations as \eqref{L0-10}, for $x \in {\mathcal R} $ and $t>\frac{1}{2}$,  we have
\begin{align}\label{0413-4}
\notag |D_{x'}^k D_{x_n} \big( D_{x_n} w_i -   {\mathcal B}_{i}^{w} \big) (x,t)|
\le &c \int_0^t \int_{0}^2 y_n^{-\beta}(\tau-\frac{1}{2})^{-\alpha} e^{-\frac{y^2_n}{t-\tau}}dy_nd\tau\\
\le &c(t-\frac{1}{2})^{\frac32  -\frac{\be}{2} -\al}.
\end{align}

$\bullet$ (case $l\geq 2$)\,\,
Direct computations show that
\begin{align}\label{0706-4}
D_{x_n} W_{mm} (x,t)   = D_{x_m}  \int_0^t \int_{\R_+} L_{nm} (x,y,t-\tau) f_2 (y,\tau) dyd\tau+{\tilde{w}}_{m}^{{\mathcal B}},
\end{align}
where
\[
{\tilde{w}}_{m}^{{\mathcal B}}:=-4 D_{x_m}^2  \int_0^t \int_{\R_+} f_2 (y,\tau) \int_{\Rn} \Ga(x' -y'-z', x_n + y_n, t-\tau)N(z', 0) dz'  dyd\tau.
\]
From \eqref{0706-3} and \eqref{0706-4}, for $l \geq 2$,     we have
\begin{align}\label{0706-5}
\notag D^l_{x_n} \big( D_{x_n} w_i -   {\mathcal B}_{i}^{w} \big) =& D^{l+1}_{x_n} V_i -\sum_{ 1 \leq m \leq n-1}  D_{x_n}^{l-2} D_{x_2}  D_{x_i} D_{x_m} W_{nm} (x,t)\\
&  + 2 D^{l-1}_{x_n} D_{x_2} D_{x_i} \Ga^* * f_2(x,t) +  D^{l-2}_{x_n}{\tilde{w}}_{m}^{{\mathcal B}}.
\end{align}
Following similar computations as in Proposition \ref{theo0525-1}, we see that for any $k \ge 0$ and $l \geq 2$
\begin{equation}\label{june16-20}
\abs{  D_{x'}^k   D^{l-2}_{x_n}  \tilde{w}^{{\mathcal B}}_i
(x,t) } \le   c \left\{\begin{array}{l}
(t-\frac{1}{2})^{\frac{3-l-\be-2\al}2 }e^{-\frac{x_n^2}{2(t-\frac{1}{2})}}, \quad l \leq 3\\
(t-\frac12)^{ -\al  }x_n^{3 -\be  -l} e^{-c_4\frac{ x_n^2}{t-\frac12} }, \quad l \geq 4.
\end{array}
\right.
\end{equation}
Its verification of \eqref{june16-20} is very similar to that of Proposition \ref{theo0525-1}, and thus the details are omitted.

Using the relations \eqref{june16-10} and  \eqref{0706-4},  from \eqref{0801-1},  \eqref{SS-lemma-10}  and the estimate of ${\mathcal B}^{w}_i$, we   have
\begin{align}\label{june16-30}
 & | D_{x_n}^{l-2}   D_{x_2}  D_{x_i} D_{x_m} W_{nm} (x,t)|\\
\le &c  (t-\frac{1}{2})^{\frac32  -\frac{\be}{2} -\al } +\left\{\begin{array}{l}
c(t-\frac{1}{2})^{1-\frac{l+\be+2\al}2 }e^{-\frac{x_n^2}{2(t-\frac{1}{2})}}, \quad l \leq 3\\
c(t-\frac12)^{ -\al  }x_n^{2 -\be  -l} e^{-c_4\frac{ x_n^2}{t-\frac12} }, \quad l \geq 4.
\end{array}
\right.
\end{align}
Combining estimates \eqref{0706-5},  \eqref{june16-20}  and \eqref{june16-30}, it follows that for $l \geq 2$,
\begin{align}\label{0413-4}
\abs{D_{x'}^k D^l_{x_n} \big( D_{x_n} w_i -   {\mathcal B}_{i}^{w} \big) (x,t)}
\leq  c \left\{\begin{array}{l}
(t-\frac{1}{2})^{1-\frac{l-2+\be+2\al}2 }e^{-\frac{x_n^2}{2(t-\frac{1}{2})}}, \quad l \leq 3\\
(t-\frac12)^{ -\al  }x_n^{2 -\be  -l} e^{-c_4\frac{ x_n^2}{t-\frac12} }, \quad l \geq 4.
\end{array}
\right.
\end{align}
We complete the proof of  the first quantity in \eqref{0525-2}.

$\bullet$ (case $i=n$)\,\,
As mentioned earlier, the tangential derivatives are rather easy to control and thus we skip its details.
Using ${\rm div  } w =0$, it is also straightforward that
\begin{align*}
\abs{D^k_{x'} D_{x_n}w_n (x,t) } =\abs{D^k_{x'} \sum_{ j=1}^{ n -1} D_{x_j} w_j (x,t) }\leq c (t-\frac12)^{\frac32-\frac\be2  -\al}.
\end{align*}
Higher derivative of $w_n$ in $x_n$ variable can be rewritten as
\[
D^{l+1}_{x_n} w_n(x,t) = -D^{l}_{x_n}\sum_{ j=1}^{ n -1} D_{x_j} w_j(x,t)
\]
\begin{equation}\label{june16-30-1}
=  -\sum_{ j=1}^{ n -1} D_{x_j}D^{l-1}_{x_n}\bke{D_{x_n} w_j(x,t)-{\mathcal B}^{w}_j  (x,t)}
 -\sum_{ j=1}^{ n -1} D_{x_j}D^{l-1}_{x_n}{\mathcal B}^{w}_j  (x,t).
\end{equation}
Using the estimates \eqref{0525-2} and Proposition \ref{theo0525-1}, we obtain via \eqref{june16-30-1}
\[
\abs{D^{l+1}_{x_n} w_n(x,t) }\le c\left\{\begin{array}{l}
(t-\frac{1}{2})^{1-\frac{l-1+\be+2\al}2 }e^{-\frac{x_n^2}{2(t-\frac{1}{2})}}, \quad l \leq 1\\
(t-\frac12)^{ -\al  }x_n^{3 -\be  -l} e^{-c_4\frac{ x_n^2}{t-\frac12} }, \quad l \geq 2.
\end{array}
\right.
\]
This completes the proof of \eqref{0616-10}.

Recalling the formulae of the pressure and its decomposition
\eqref{formulas-p-1}-\eqref{badterm-pi}
and using  \eqref{0515-1}, we get
\begin{align*}
|  D_{x}^k\Pi^{{\mathcal G}}(x,t)| & \le   c_n \int_0^t \int_{{\mathbb R}^n_+}
     f_2(y, \tau) (t-\tau)^{-\frac12}     e^{-\frac{y_n^2}{t-\tau}}
     |D_{x}^k D_{x_2}D_{x_n}  N( x'-y',x_n)| dyd\tau\\
 & \quad   +  c_n \int_0^t \int_{{\mathbb R}^n_+}
     f_2(y, \tau) (t-\tau)^{-\frac12}     e^{-\frac{y_n^2}{t-\tau}}|J^2_1(x'-y', t-\tau)|  dyd\tau\\
&  : = \Phi^{{\mathcal G}}_1 (x,t) + \Phi^{{\mathcal G}}_2(x,t).
\end{align*}

Since $|x' -y'| \geq 1$ for $y' \in B'_1$ and $|x_2| \geq 2$, from \eqref{0515-1}, we have
\begin{align*}
|\Phi^{{\mathcal G}}_{1}(x,t)| &  \leq  c  \int_0^t  (\tau-\frac{1}{2})^{-\al}    \int_0^1 y_n^{-\be}\int_{|y'| \leq 1}
     g^{{\mathcal T}}(y') (t-\tau)^{-\frac12}     e^{-\frac{y_n^2}{t-\tau}}  |x'-y'|^{-n-k} dy' dy_n d \tau \\
&  \leq  c  \int_0^t  (\tau-\frac{1}{2})^{-\al}  (t -\tau)^{-\frac12}   \int_0^1 y_n^{-\be} e^{-\frac{y_n^2}{t-\tau}}  dy_n d \tau
\leq c (t-\frac12)^{1 -\frac{\be}2 -\al}
\end{align*}
and
\begin{align*}
|\Phi^{{\mathcal G}}_{2}(x,t)| &  \leq  c  \int_0^t (\tau-\frac{1}{2})^{-\al}   \int_0^1 y_n^{-\be}\int_{|y'| \leq 1}
     g^{{\mathcal T}}(y')     e^{-\frac{y_n^2}{t-\tau}}   dy' dy_n d \tau
\leq c  (t-\frac12)^{\frac32 -\frac{\be}2 -\al}.
\end{align*}
Thus, we deduce \eqref{0525-1}, and
 the proof Proposition \ref{prop0525-3} is completed.

\begin{rem}
 We remark that from the same estimate in \eqref{L0-10}, for $x \in {\mathcal R} $ and $t>\frac{1}{2}$, we obtain
 the following estimate:
 \begin{align}\label{0709-1}
\abs{w(x,t) } \leq c (t-\frac{1}{2})^{\frac32  -\frac{\be}{2} -\al }.
\end{align}
Since its verification is rather straightforward, the details are omitted.
\end{rem}

\subsection{Proof of Proposition \ref{prop0525-3-1}}
\label{proofprop0525-3-1}

From the equations for $i<n$, it follows that
\begin{align*}
D_t w_i - D_{x_n} w_i^{{\mathcal B}} + D_{x_i} \Pi^{{\mathcal B}} = f_i  +\De' w_i + D_{x_n}  w_i^{{\mathcal G}} - D_{x_i} \Pi^{{\mathcal G}}.
\end{align*}
Using Proposition \ref{prop0525-3} and Lemma \ref{theo0521-1}, we obtain
\eqref{May27-90}. On the other hand, the equation of $w_n$ can be rewritten as
\[
D_t w_n- D^2_{x_n} w_n+ D_{x_n} \Pi^{{\mathcal B}}  = f_n  +\De' w_n - D_{x_n} \Pi^{{\mathcal G}},
\]
which yields, due to Proposition \ref{prop0525-3} and divergence free condition, the estimate \eqref{0617-100}.
This completes the proof of Proposition \ref{prop0525-3-1}.

\section{Proofs of Theorems for Stokes system}\label{proofofmaintheo}
\setcounter{equation}{0}

\subsection{Proof of Theorem \ref{maintheo0503}}

\label{prooftheomaintheo0503}

Assume that  $0 < \al<1$ and $0 <\be<\frac12$. Let $\frac1{q_1} = 1 -\de$ and $\frac1{p_1} = \frac12 +\ep$ for sufficiently small $ 0 < \frac{n\ep}2 < \de$ such that   $  \frac12   <\frac{1}{p_1}$,  $ \max(\frac12, \al) < \frac1{q_1}$  and
 \begin{align}\label{0504-8}
   \frac2{q_1} + \frac{n}{ p_1} = 2-  2\de +\frac{n}2+ n\ep <  \frac{n}2 +2.
 \end{align}
Note that $f_2 \in L^{q_1} (0,\infty; L^{p_1} (\R_+))$. Then, from \eqref{maximal-SS}, \eqref{loworder-SS} and \eqref{mixednorm-f}, we get
\begin{align}\label{energy-bound-3}
\| w\|_{L^\infty(0,1; L^2(\R_+))} < \infty,
\end{align}
\begin{equation}\label{energy-bound-2}
 \|D_{x} w\|_{L^2 (\R_+ \times (0, 1))}   < \infty.
\end{equation}
From \eqref{energy-bound-2} and \eqref{energy-bound-3}, we obtain \eqref{energy-bound}.
Since $ q> 6$ with $ 1 +\frac3{2q} < \al +\frac{\be}2$, using Proposition \ref{prop0525-3}  for $i \neq 2$ and  Proposition \ref{theo0521-1}, we obtain
\begin{align*}
&\min \bket{\| D_{x_n} w_i\|^q_{L^q_t( 0,1; L^q_x( A_{i} \times(0, 1))}, \,\, \| D_{x_n} w_2\|^q_{L^q_t( 0,1; L^q_x( B_{2} \times(0, 1))}}\\
\geq & c  \int_{\frac{1}{2}}^1  (t-\frac12)^{(1 -\frac{\be}2 -\al)q + \frac{1}{2}}  dt - \int_{\frac{1}{2}}^1  (t-\frac12)^{(\frac32 -\frac{\be}2 -\al)q} dt = \infty.
\end{align*}
Hence, we obtain   \eqref{0504-2}. We complete the proof of Theorem \ref{maintheo0503}.

\subsection{Proof of Theorem \ref{maintheo0503-2}}

Let $\frac2q  +\frac{n}{p} + 1 > 2\al +n \be$ for $ \frac{n}{n-1} < p$ and $\frac1{\be} < \frac{n-1}{n} p$. We take $ q_1 < \frac1{\al}$ and $p_1< \frac1{\be}$ satisfying $ \frac2q +\frac{n}{p} +1> \frac2{q_1}  + \frac{n}{p_1}  $, and then by Proposition \ref{prop0503}, it follows that
\begin{equation*}
\| \Pi\|_{L^{q}(0, 1; L^p (\R_+))} < c \| f_2\|_{L^{q_1} (0, 1; L^{p_1} (\R_+))}.
\end{equation*}
Hence, we obtain \eqref{pi-bound}.
We note from \eqref{0525-1} an \eqref{june17-10} that there exists $\delta\in (\frac{1}{2}, 1)$ such that  if $t\in (\frac{1}{2}, \delta)$, then
\[
\abs{\Pi(x,t)}\ge c (t-\frac{1}{2})^{\frac{1}{2}-\frac{\be}{2}-\al}(1-c(t-\frac{1}{2})^{\frac{1}{2}})\ge c(t-\frac{1}{2})^{\frac{1}{2}-\frac{\be}{2}-\al}.
\]
Therefore, we obtain
\begin{align*}
\| \Pi\|^q_{L^q (\{|x'| >2 \} \times (a ,b) \times (0,1))}
\geq c\int_{\frac12}^{\delta} (t -\frac12)^{(\frac12 -\frac{\be}2 -\al)q} dt
=\infty
\end{align*}
for $\frac12 +\frac1q < \al +\frac{\be}2$. Hence, we obtain \eqref{0504-4}. We complete the proof of Theorem \ref{maintheo0503-2}.
\qed

\subsection{Proof of Theorem \ref{theo2}}

Let  $0 < \alpha, \, \be < 1 $ and $w$ be a  solution of \eqref{maineq} defined by \eqref{int-0909-1}.
We recall  the decomposition of $w$ in \eqref{V-eqn} and \eqref{W-eqn}, i.e. $w = V + W$.  We assume that $\ep_0 = 3 -\be -2\al\in  (0, 2) $.
It is clear that
\[
V \in C^{\infty}_{x,t} (  {\mathcal R} \times (0,1)).
\]
Applying the proof of Proposition \ref{theo0521-1},  Proposition \ref{theo0525-1} and Proposition \ref{prop0525-3}, for $x \in {\mathcal R} $ and $t>\frac{1}{2}$, if $\ep_0 \in (0, 1)$, then we have
\begin{align*}
x_n^{1 -\ep_0 } |D_{x} W(x,t)| & \leq   c x_n^{1 -\ep_0}(t-\frac12)^{1 -\frac{\be}2 -\al} e^{-\frac{x_n^2}{t-\frac12}} + cx_n^{1 -\ep_0 }\\
&\le c(\frac{x_n^2}{t-\frac12})^{\al+\frac{\be}2-1} e^{-\frac{x_n^2}{t-\frac12}} + cx_n^{2\al+\be-2  }\le c,
\end{align*}
and  if $\ep_0 \in (1, 2)$, then we have
\begin{align*}
x_n^{2 -\ep_0 } |D^2_{x} W(x,t)| & \leq   c x_n^{2 -\ep_0}(t-\frac12)^{\frac12 -\frac{\be}2 -\al} e^{-\frac{x_n^2}{t-\frac12}} + cx_n^{2 -\ep_0 }\\
&\le c(\frac{x_n^2}{t-\frac12})^{\al+\frac{\be}2-\frac12} e^{-\frac{x_n^2}{t-\frac12}} + cx_n^{2\al+\be-2  }\le c,
\end{align*}
which implies that
\begin{align}\label{holderspace}
W \in L^\infty (0, 1; C^{\ep_0} (  {\mathcal R}  \times  (0,1)))
\end{align}
(see the proof of Theorem 4.1 in \cite{JK}).

For \hoe continuity for temporal variable, we set $s,t$ with $t > s > \frac12$.
\begin{align*}
W_i(x,t) - W_i(x,s) &  = \int_\frac12^s \int_{\R_+} \big( L_{i2} (x,y, t-\tau) - L_{i2}(x,y,s-\tau) \big) f_2(y,\tau) dyd\tau\\
& \quad + \int_s^t \int_{\R_+}  L_{i2} (x,y, t -\tau) f_2(y,\tau) dyd\tau\\
& = I_1 +I_2.
\end{align*}
We firstly estimate $I_2$.
\begin{align}\label{est-j2}
\nonumber \abs{I_2}&\le \int_s^t \int_{\R_+}  \frac{ e^{-\frac{y_n^2}{t- \tau}}}{ ( |x -y|^2 +  t - \tau)^{\frac{n}2}  }    y_n^{-\be} (\tau-\frac{1}{2})^{-\al}  dy d\tau\\
\nonumber &\leq  \int_s^t \int_{\R_+}    e^{-\frac{y_n^2}{t- \tau}}    y_n^{-\be} (\tau-\frac{1}{2})^{-\al}  dy d\tau\\
&\leq  \int_s^t (t -\tau)^{ \frac12 -\frac{\be}2} (\tau-\frac{1}{2})^{-\al}  d\tau
\leq c ( t-s)^{\frac32 -\frac{\be}2 -\al}.
\end{align}

Next, we estimate $I_1$.
\[
I_1 = \int_\frac12^s \int_{\R_+} \int_0^1 D_t L_{i2} (x,y, \te t + (1 -\te)s- \tau) ( t -s)  f_2(y,\tau) d\te dyd\tau
\]
\[
 \leq c(t-s) \int_\frac12^s \int_{\R_+} \int_0^1  \frac{ e^{-\frac{y_n^2}{\te t + (1 -\te)s- \tau}}f_2(y,\tau)}{( \te t + (1 -\te)s- \tau) ( |x -y|^2 + \te t + (1 -\te)s- \tau)^{\frac{n}2}  }     d\te  dyd\tau
\]
\[
\leq c(t-s) \int_\frac12^s \int_0^1 \int_0^1  \frac{ e^{-\frac{y_n^2}{\te t + (1 -\te)s- \tau}}}{( \te t + (1 -\te)s- \tau) }  y_n^{-\be} (\tau-\frac{1}{2})^{-\al} d\te  dy_n d\tau
\]
\[
 \leq c(t-s) \int_\frac12^s   \int_0^1   ( \te t + (1 -\te)s- \tau)^{-\frac12 -\frac{\be}2  }     (\tau-\frac{1}{2})^{-\al} d \te    d\tau
\]
\[
= \int_\frac12^s (\tau-\frac{1}{2})^{-\al} (s -\tau)^{\frac12-\frac{\be}2} \int_0^{\frac{t-s}{s -\tau}} (\te + 1)^{ -\frac12 -\frac{\be}2}d\te d\tau
\]
\[
\le \int^{2s -t}_{\frac{1}{2}} \int_\frac12^{\frac{t-s}{s -\tau}}\cdots d\te d\tau+\int_{2s -t}^s \int_0^{\frac{t-s}{s -\tau}}\cdots d\te d\tau :=I_{11}+I_{12}.
\]
Using $\frac{t-s}{s -\tau}\le 1$ for $\tau<2s-t$, it follows for $I_{11}$ that
\begin{align}\label{est-j11}
\nonumber I_{11}&\le  c (t-s)    \int_\frac12^{2s -t} (\tau-\frac{1}{2})^{-\al} (s -\tau)^{-\frac12-\frac{\be}2}  d\tau \\
&\le c (t-s)^{\frac12-\frac{\be}2}    \int_\frac12^{2s -t} (\tau-\frac{1}{2})^{-\al}  d\tau\le c (t-s)^{\frac32-\frac{\be}2-\al}.
\end{align}
On the other hand, since $ \int_0^{\frac{t-s}{s -\tau}} (\te + 1)^{ -\frac12 -\frac{\be}2}d\te\le c(\frac{t-\tau}{s -\tau})^{\frac12 -\frac{\be}2}$, we have
\begin{equation}\label{est-j12}
I_{12}\le  c (t -s)^{\frac12 -\frac{\be}2} \int_{2s -t}^s (\tau-\frac{1}{2})^{-\al}
\le c (t-s)^{\frac32-\frac{\be}2-\al}.
\end{equation}
Adding up estimates \eqref{est-j11} and \eqref{est-j12}, we obtain
\begin{align}\label{est-j1}
I_1 \leq c    (t-s)^{\frac12\ep_0}.
\end{align}
Hence, it follows from \eqref{est-j1} and \eqref{est-j2} that
\begin{align}\label{holdertime}
W \in L^\infty ( {\mathcal R}; C^{\frac12\ep_0 } ( 0,1 )).
\end{align}
From \eqref{holderspace} and \eqref{holdertime}, we obtain \eqref{0830-1}.

Next, we prove \eqref{0830-2}.
We note first due to \eqref{0515-1} that
\begin{align*}
&L_{i2}(x,y,t-s) = D_{x_i} \int_0^{x_n} \int_{\Rn}\Ga(x -y^* -z, t-s) D_{z_2} N(z) dz\\
&\quad =  c (t-s)^{-\frac12}   \int_0^{x_n} e^{-\frac{(x_n + y_n - z_n)^2}{t-s}} \big( D_{x_i} D_{x_2}  N( x'-y',z_n)    + J_{20}(x'-y', t-s) \big) dz_n.
\end{align*}
For $x' \in A_{i}$ and $|y'| < 1$, we estimate the second term in the above as follows:
\begin{align*}
(t-s)^{-\frac12}    \abs{\int_0^{x_n} e^{-\frac{(x_n +y_n - z_n)^2}{t-s}}   J_{20}(x'-y', t-s)  dz_n}
 \leq c   \int_0^{x_n} e^{-\frac{( y_n + z_n)^2}{t-s}} dz_n
  \leq c (t -s)^\frac12e^{-\frac{y_n^2}{t-s}}.
\end{align*}
On the other hand, the first term has a lower bound. Indeed, since  $  D_{x_i} D_{x_2}  N( x'-y',z_n)  \geq \phi_i(x', x_n)$ for $x' \in A_{i}$, $|y'| \leq 1$ and $0 < z_n \leq x_n$, the first term is lower bounded by
\begin{align*}
&(t-s)^{-\frac12}   \int_0^{x_n} e^{-\frac{(x_n + y_n-z_n)^2}{t-s}}  \phi_i(x,z_n)    dz_n\\
  \geq &c    (t-s)^{-\frac12}   \int_0^{x_n} e^{-\frac{( y_n+z_n)^2}{t-s}}     dz_n  \phi_i(x',x_n) \\
  \geq  &c  e^{-\frac{ y_n^2}{t-s}} \int_0^{\frac{x_n}{\sqrt{t-s}}} e^{-z_n^2} dz_n \phi_i(x',x_n) .
\end{align*}
Recalling $W_i$ in \eqref{W-eqn}. Since $W_i(x', 0,t) =0$,  we have  for $x' \in A_{i}$
\begin{align*}
|W_i(x,t) - W_i(x',0,t)| & \geq c  \int_{\frac12}^t   \int_0^1 e^{-\frac{ y_n^2}{t-s}}   y_n^{-\be} (s-\frac12)^{-\al} dy_nds |\phi_i(x,x_n) |\\
& \quad - c\int_{\frac12}^t \int_0^1 ( t-s)^\frac12  e^{-\frac{y_n^2}{t-s}} y_n^{-\be} (s-\frac12)^{-\al} dy_nds\\
& \geq c(t-\frac12)^{\frac32 -\frac{\be}2 -\al} - c(t-\frac12)^{2 -\frac{\be}2 -\al}.
\end{align*}
Therefore, for $x' \in A_i$ and  $x_n^2 < t -\frac12 < 4 x_n^2$, we obtain
\begin{align*}
|W_i(x,t) - W_i(x',0,t)| \geq  c (x_n^{3 -\be -2\al} -  x_n^{4 -\be -2\al}),
\end{align*}
which deduce that
\begin{align}\label{0821-1}
W_i \notin L^\infty_t C^{ \ep }_x ( A_i \times  (0, 1)  \times (0, 1))
\end{align}
for   $  \ep > \ep_0 = 3 -\be -2 \al$ if $\ep_0 \in (0, 1)$.

Let $\ep_0 \in (1,2)$. From the proof of  Proposition \ref{prop0525-3},  for $x \in {\mathcal R} $ and $t>\frac{1}{2}$, we have
\begin{align*}
x_n^{2-\ep} | D_{x} D_{x_2} W_{i}^{{\mathcal G}} (x,t) \big)|
& \le  c  x_n^{2-\ep}
(t-\frac{1}{2})^{1-\frac{\be}2 -\al}e^{-\frac{x_n^2}{8(t-\frac{1}{2})}}\\
& \le  c  x_n^{4-\ep -\be -2\al}\\
& < \infty \quad \mbox{for} \quad \ep_0 < \ep < 4-\be -2\al <2.
\end{align*}
This implies $ D_{x_2} W_{i}^{{\mathcal G}} \in L^\infty(0, 1; C^{\ep} (  {\mathcal R} ))$.

From Mean-value Theorem and  Proposition  \eqref{theo0525-1}, for $\sqrt{t -\frac12} \leq c_l x_n$, we have
\begin{align*}
|{\mathcal B}^{w}_i  (x', x_n,t) - {\mathcal B}^{w}_i  (x', \frac12 x_n,t)| & = \frac12 x_n | D_{x_n} {\mathcal B}^{w}_i  (x', \xi_n,t)| \quad \mbox{for some} \quad \xi_n \in [\frac12 x_n, x_n]\\
&\ge
c  \frac12 x_n (t -\frac12)^{\frac12 -\frac{\be}2 -\al} e^{-\frac{x_n^2}{2(t-\frac12)}}.
\end{align*}
Therefore, for $|x'| \geq 2$ and  $c'_l x_n < \sqrt{ t -\frac12 }< c_l x_n$, we obtain
\begin{align*}
|{\mathcal B}^{w}_i  (x', x_n,t) - {\mathcal B}^{w}_i  (x', \frac12 x_n,t)|&\ge
c   x_n^{2 -\be -2\al}.
\end{align*}
This implies that $ D_{x_n} {\mathcal B}^{w}_i \notin L^\infty(0,1; C^{\ep -1}(  {\mathcal R}) )$ for $\ep > \ep_0$. Since $D_{x_n} W_i =  D_{x_2} W_{i}^{{\mathcal G}} + {\mathcal B}^{w}_i$, We have that   $ W_i \notin L^\infty(0,1; C^{\ep }(  {\mathcal R}))$ for $\ep > \ep_0 $ if $\ep_0 \in (1,2)$.

Next, we will show that $W_i \notin  C^{ \frac{\ep}{2} }_t L^\infty_x( A_i \times  (0, 1)  \times (0, 1))$. Indeed, for $t > \frac12$ we consider
\begin{align*}
& D_{x_i}D_{x_2}\int_0^{x_n} \int_{\Rn} \Ga(x-y^* -z,t) N(z) dz\\
 =& \int_0^{x_n} \Ga_1(x_n + y_n -z_n, t) D_{z_i} D_{z_2} N(z) dz
 + \int_0^{x_n} \Ga_1(x_n + y_n -z_n, t) J_{20}(z,t) dz\\
 =& I_1 + I_2,
\end{align*}
where we used Lemma \ref{lemma0709-1}.
We first estimate $I_2$ as follows:
\begin{align*}
| I_2|  & \leq  c   \int_0^{x_n} e^{-\frac{(x_n + y_n -z_n)^2}t}  dz_n  \leq  c  e^{ -\frac{y_n^2} t} \int_0^{x_n} e^{-\frac{z_n^2}t}  dz_n \leq  c t^\frac12 e^{ -\frac{y_n^2} t}.
\end{align*}
Secondly,
for $x_n^2 \geq t -\frac12$, we have a lower bound for $I_1$.
\begin{align*}
I_1 & \geq  c   t^{-\frac12} \int_0^{x_n} e^{-\frac{(x_n + y_n -z_n)^2}t}  dz_n
 \geq  c    e^{-\frac{y_n^2}t} \int_0^{\frac{x_n}{\sqrt{t}}} e^{- z_n^2 }  dz_n \geq c e^{-\frac{y_n^2}t}.
\end{align*}
Collecting above estimates and noting that $ W_i(x,\frac12) =0$, we obtain for $x_n^2 \geq t -\frac12$
\begin{align}\label{june23-10}
\abs{W_i(x,t) }=\nonumber &\abs{W_i(x,t) - W_i(x,\frac12)}  = \abs{\int_{\frac12}^t L_{i2}(x, y, t-\tau) y_n^{-\be} \tau^{-\al} dyd\tau}\\
\nonumber  \geq & c \int_\frac12^t \int_0^1 y_n^{-\be} s^{-\al} e^{-\frac{y_n^2}{(t-\tau)}} dy_n d\tau - c \int_\frac12^t \int_0^1 y_n^{-\be} \tau^{-\al} (t -\tau)^\frac12  e^{-\frac{y_n^2}{(t-\tau)}} dy_n d\tau\\
\nonumber  \geq & c \int_\frac12^t \tau^{-\al} (t-\tau)^{\frac12 -\frac{\be}2}   d\tau - c \int_\frac12^t \tau^{-\al} (t-\tau)^{1 -\frac{\be}2}   d\tau\\
\geq & c \big(  (t-\frac12)^{\frac32 -\frac{\be}2 -\al}    - c(t-\frac12)^{2 -\frac{\be}2 -\al} \big),
\end{align}
where we used \eqref{SS-lemma-10}.
The estimate \eqref{june23-10} implies that
\begin{align}\label{0821-2}
W_i \notin L^\infty(A_i \times (0, 1); \dot C^{\frac12 \ep } (0, 1))
\end{align}
for $\ep > \ep_0$. From \eqref{0821-1} and \eqref{0821-2}, we obtain \eqref{0830-2}. We complete the proof of Theorem \ref{theo2}.

%%%%%%%%%%%%%%%%%%%%%%%%%%%%%%%%%%%%%%%%%%%%
%%%%%%%%%%%%%%%%%%%%%%%%%%%%%%%%%%%%%%%%%%%%
%%%%%%%%%%%%%%%%%%%%%%%%%%%%%%%%%%%%%%%%%%%%

\section{Singular solutions for Navier-Stokes equations }
\label{proof-NSE}
\setcounter{equation}{0}

In this section we provide the proof of Theorem \ref{thm-NSE}.

We look for a solution of the Navier-Stokes equations \eqref{NSE-eq} of the form $u=w+v$ and $p = \Pi + q$, where $(w,\Pi)$ is the singular solution of the Stokes system \eqref{maineq}
satisfying \eqref{Sept08-55}. Thus, it suffices to establish existence of solution that is  less singular than $(w, \Pi)$ for the following perturbed Navier-Stokes equations in
$\R_+\times (0,1)$:
\begin{equation}\label{CCK-Feb7-10}
v_t-\Delta v+\nabla q+{\rm div}\,\left(v\otimes v+v\otimes
w+w\otimes v\right)=-{\rm div}\,(w\otimes w), \quad {\rm
div} \, v =0
\end{equation}
with homogeneous initial and boundary data, i.e.
\begin{equation}\label{pnse-bdata-20}
v(x,0)=0,\qquad v(x,t)=0 \,\,\mbox{on} \,\,\{x_n=0\}.
\end{equation}

The first step is to construct a solution of \eqref{CCK-Feb7-10}-\eqref{pnse-bdata-20} in the class of functions $L^r(\R_+ \times (0, 1)) \cap L^\infty (0, 1; L^2 (\R_+))$, where $r$ is the number imposed in \eqref{Sept08-10}. In order to do that, we use the following proposition (see \cite[Theorem 1.2]{CJ}).
\begin{prop}
\label{prop-stokes}
Let    $1<p < \infty, \, 1 <  q< \infty$,   $ p \geq p_1$ and $q\geq q_1$.
If  $F\in L^{q_1}(0,1; L^{p_1}(\R_+))
$, $(\frac{n}{2p_1}-\frac{n}{2p}) + \frac{1}{q_1} -\frac{1}{q}  \leq  \frac12$  with $F|_{x_n =0} =0$. Then, there is a unique weak solution $u\in L^q(0,1;L^p(\R_+))$ to the Stokes equation \eqref{maineq} with  $ f = {\rm div   } F$, which  satisfies the estimate
\begin{align}\label{0919-3}
\| v\|_{L^q(0,1;L^p(\R_+))}\leq c\|F\|_{L^{q_1}(0,1; L^{p_1}(\R_+))}.
\end{align}
\end{prop}

Now, we
adopt an iterative scheme for \eqref{CCK-Feb7-10}, which is
formulated as follows: For a positive integer $m\ge 1$
\begin{align*}
&v^{m+1}_t-\Delta v^{m+1}+\nabla q^{m+1}=-{\rm
div}\,\left(v^{m}\otimes v^{m}+v^{m}\otimes w+w\otimes
v^{m}+w\otimes w \right),\\
& \qquad \qquad \qquad \qquad \qquad {\rm div} \, v^{m+1} =0
\end{align*}
with conditions \eqref{pnse-bdata-20} i.e.
\[
v^{m+1}(x,0)=0\qquad
v^{m+1}(x,t)=0\,\,\mbox{on} \,\,\{x_n=0\}.
\]
We set $v^1=0$.

For convenience, we denote that $L^q(0, 1; L^p (\R_+)) = L^q_tL^p_x$. Furthermore, if $p=q$, then we abbreviate $L^q_tL^q_x$ as $L^q$, unless any confusion is to be expected.
Due to Proposition \ref{prop-stokes}, we have
\begin{align}\label{est-v2-10}
& \|   v^2\|_{L^2}  \leq c     \||w|^2\|_{L^{2} }
  \leq c   \|w\|^2_{L^4   }.%,\\
& \|   v^2\|_{L^r  }  \leq c     \||w|^2\|_{L^{\frac{(n+2)r}{ n+2 + r}} }
  \leq c   \|w\|^2_{L^{\frac{2(n+2)r}{ n+2 + r}}   },
\end{align}
\begin{equation}\label{0915-1}
\|   v^{m+1}\|_{L^r  }  \leq c    \big( \|w \otimes w\|_{L^{\frac{(n+2)r}{ n+2 + r}} } + \|v^m \otimes w\|_{L^{\frac{(n+2)r}{ n+2 + r}} } + \|v^m \otimes v^m\|_{L^{\frac{(n+2)r}{ n+2 + r}} } \big),
\end{equation}
\begin{equation}\label{0915-2}
\|v^m\|_{L^2  }  \leq c    \big(   \|w \otimes w\|_{L^2 } + \|v^m \otimes w\|_{L^2  } + \|v^m \otimes v^m\|_{L^2 } \big).
\end{equation}
Since $n+2 < r$, we take $0 < \eta_1 =\frac{nr}{(r-2)(n+2)} \in (0, 1)$.
Then, we have
\begin{align*}
\|v^m \otimes v^m\|_{L^{\frac{(n+2)r}{ n+2 + r}} }\leq \|v^m  \|_{L^r } \|v^m  \|_{L^{n+2} }
 \leq \|v^m  \|^{1 +\eta_1}_{L^{r} } \|v^m  \|^{1 -\eta_1}_{L^2 }  \leq \|v^m  \|^{1 +\eta_1}_{L^{r} } \|v^m  \|^{1 -\eta_1}_{L^2 }.
\end{align*}
Let $r_2$ be the number with $\frac14 = \frac1r +\frac1{r_2}$. Since $r>8$, it follows that $2<r_2<r$ and thus
\begin{align*}
\|v^m \otimes v^m\|_{L^2 }\leq \|v^m  \|^2_{L^4 }
 \leq \|v^m  \|^{2(1 -\eta_2)}_{L^{r} } \|v^m  \|^{2\eta_2}_{L^2 }  \leq \|v^m  \|^{2(1 -\eta_2)}_{L^{r} } \|v^m  \|^{2\eta_2}_{L^2},
\end{align*}
where $\eta_2=\frac{r(r_2-2)}{r_2(r-2)}$. We note that $\eta_2\ge \frac{1}{2}$.
Using computations above, \eqref{0915-1} and \eqref{0915-2} are controlled as follows:
\begin{align}
\label{0904-1}& \|   v^{m+1}\|_{L^r  }  \leq c   \big(  \|w  \|^{1 +\eta_1}_{L^{r} } \|w  \|^{1 -\eta_1}_{L^2  } + \|v^m  \|^{1 +\eta_1}_{L^{r} } \|w  \|^{1 -\eta_1}_{L^2} + \|v^m  \|^{1 +\eta_1}_{L^{r} } \|v_m  \|^{1 -\eta_1}_{L^2 }  \big),\\
\label{0904-2}& \|v^m\|_{L^2 }
  \leq c      \big( \|w  \|^2_{L^4 } +\|v^m  \|^{2(1 -\eta_2)}_{L^{r} } \|v^m  \|^{2\eta_2}_{L^2 } + \|w  \|^{2(1 -\eta_2)}_{L^{r} } \|v^m  \|^{2\eta_2}_{L^2 } \big).
\end{align}
By \eqref{maximal-SS}, we have  $ A :=    \|w  \|_{L^{r} } + \|w  \|_{L^2} \leq ca$, where $a>0$ is defined in {\bf Assumption}  \ref{force-f}.

Taking $a >0$ small such that
$A<\frac{1}{4c}$, where $c$ is the constant in
\eqref{est-v2-10}-\eqref{0904-2} such that
\begin{align*}
   \| v^2\|_{L^2} + \|  v^2\|_{L^r }  < A.
\end{align*}
Moreover, iterative arguments show for any $m$ that
\begin{align}\label{0531-1}
\notag&    \| v^{m+1}\|_{L^2} + \|  v^{m+1}\|_{L^r} \\
 \notag & \quad \leq c \big(\|w  \|^{2(1 -\eta_2)}_{L^{r} } \|w  \|^{2\eta_2}_{L^2} +\|v^m  \|^{2(1 -\eta_2)}_{L^{r} } \|v^m  \|^{2\eta_2}_{L^2} + \|w  \|^{2(1 -\eta_2)}_{L^{r} } \|v^m  \|^{2\eta_2}_{L^2 } \big) \\
 & \quad   \leq 4c A^2 <A.
\end{align}

Next, we will show that $v^m$ converges in $L^2\cap L^r$.
For simplicity, we denote $V^{m+1}:=v^{m+1}-v^{m}$ and $Q^{m+1}:=q^{m+1}-q^{m}$ for $m\ge
1$. We then see that $(V^{m+1}, Q^{m+1})$ solves
\[
V^{m+1}_t-\Delta V^{m+1}+\nabla Q^{m+1}=-{\rm
div}\,\left(V^{m}\otimes v^{m}+v^{m-1}\otimes V^{m}+V^{m}\otimes
w+w\otimes V^{m}\right),
\]
\[
{\rm div} \, V^{m+1} =0
\]
with homogeneous initial and boundary data, i.e. $V^{m+1}(x,0)=0$
and $V^{m+1}(x,t)=0$ on $\{x_n=0\}$. Taking sufficiently small
$a>0$ such that $A < \frac1{6c}$.
Since $\eta_2  \geq \frac12$, it follows  from \eqref{0904-1}, \eqref{0904-2} and  \eqref{0531-1} that
\begin{align}\label{0920-1}
\begin{split}
  \| V^{m+1}\|_{L^2} + \|  V^{m+1}\|_{L^r }
&\leq c \big(\|v^m  \|^{2(1 -\eta_2)}_{L^{r} } \|V^m  \|^{2\eta_2}_{L^2 } + \|w \|^{2(1 -\eta_2)}_{L^{r} } \|V^m \|^{2\eta_2}_{L^2} \big)\\
&\leq 3c A^{2(1 -\eta_2)} \|V^m  \|^{2 \eta_2}_{L^2}\\
&\leq 3c A  \|V^m  \|_{L^2 }
<\frac12 \norm{V^m}_{L^2 }.
\end{split}
\end{align}

Hence, we obtain
\begin{align*}
  \| V^{m+1}\|_{L^2} + \|  V^{m+1}\|_{L^r }
<\frac12 \big( \norm{V^m}_{L^2 } +\|  V^{m}\|_{L^r } \big).
\end{align*}

Therefore, there exists $v \in L^r  \cap L^2 $ such that $v^m $ converges to $v $ in $L^r  \cap L^2 $, which $v$
solves \eqref{CCK-Feb7-10}-\eqref{pnse-bdata-20}
in the sense of distributions with corresponding  pressure $q$.

For the uniqueness, we assume that $v_1$ is another solution of  \eqref{CCK-Feb7-10}  with initial-boundary condition \eqref{pnse-bdata-20} such that $\| v_1\|_{L^r} + \| v_1\|_{L^2} < A$. Let   we denote $V:=v_1-v$ and $Q:=q_1-q$. We then see that $(V, Q)$ solves
\[
V_t-\Delta V+\nabla Q=-{\rm
div}\,\left(V\otimes v_1+v_1\otimes V+V\otimes
w+w\otimes V\right),\qquad
{\rm div} \, V =0
\]
with homogeneous initial and boundary data, i.e. $V(x,0)=0$
and $V(x,t)=0$ on $\{x_n=0\}$. With same estimate to \eqref{0920-1}, we obtain
\begin{align*}
  \| V\|_{L^2} + \|  V\|_{L^r }
<\frac12 \big( \norm{V}_{L^2 } +\|  V\|_{L^r } \big).
\end{align*}
Hence, we obtain $ V \equiv 0$.

By uniqueness,  $v$ is represented by
\begin{align}\label{0920-3}
v(x,t) = \int_0^t \int_{\R_+} K(x,y, t-s) {\mathbb P} {\rm div} \, F(y,s)dyds,
\end{align}
where $F = \left(v\otimes v+v\otimes
w+v\otimes w+w\otimes w\right)$ and  $K$ is introduced in \eqref{formulas-k}.
Since $\rabs{F}_{x_n =0} =0$, it is known that there exists a tensor $\calF $ such that ${\mathbb P} {\rm div} \, F = \na \cdot \calF$  with ${\mathcal F}_{in}|_{x_n =0} =0, \, 1 \leq i \leq n$ and $\| \calF\|_{L^p(\R_+)} \leq c \| F\|_{L^p (\R_+)}, \, 1 < p < \infty$. (see  \cite{KS2002}  and \cite{S03}). Hence, we have
\begin{align*}
v(x,t) = -\int_0^t \int_{\R_+} \na K(x,y, t-s) {\mathcal  F}(y,s)dyds
\end{align*}
and from \eqref{0801-1}, we have
\begin{align*}
|v(x,t)| & \leq c\int_0^t (t -s)^{-\frac12 - \frac{n}{r}} \| F(s) \|_{L^{\frac{r}2} (\R_+)} ds\\
& \leq c\int_0^t (t -s)^{-\frac12 - \frac{n}{r}} \| F(s) \|_{L^{\frac{r}2} (\R_+)} ds\\
& \leq c t^{1 -\frac{n+2}{r}} \| F\|_{L^{\frac{r}2}}.
\end{align*}
Since $ r > n+2$, we get $\| v\|_{L^\infty} < \infty$. Using the complex interpolation, we get $\| v\|_{L^p} < \infty$ for $2 \leq p \leq \infty$.
This implies that
\begin{align}\label{0920-4}
\begin{split}
 & \| \na  v\|_{L^{\frac{r}2}  }  \leq c     \|w \otimes w\|_{L^{\frac{r}2} } +  \|v \otimes w\|_{L^{\frac{r}2} } +  \|v \otimes v\|_{L^{\frac{r}2} } < \infty,\\
& \|\na  v\|_{L^2 }    \leq c       \|w \otimes w  \|^2_{L^2 } +\|w \otimes  v \|^2_{L^2 } + \|v \otimes v  \|^2_{L^2 } < \infty.
\end{split}
\end{align}
From representation \eqref{0920-3} and estimate \eqref{0801-1}, we have
\begin{align*}
\|v(t)\|_{L^2(\R_+)} & \leq c\int_0^t  \| {\mathbb P}  {\rm div  } F(s) \|_{L^{ 2} (\R_+)} ds\\
& \leq c\int_0^t \big( \| (\na w)w\|_{L^2 (\R_+)}  + \| (\na v)w\|_{L^2 (\R_+)} + \| (\na w)v\|_{L^2 (\R_+)} + \| (\na v)v\|_{L^2 (\R_+)} \big)  ds\\
& \leq c  \big( \| \na w\|_{L^2_tL^4_x  }\|  w\|_{L^2 L^4_x (\R_+)}  + \| \na w\|_{L^2_tL^4_x  }\|  v\|_{L^2 L^4_x (\R_+)}\\
  & \qquad + \| \na v\|_{L^2_tL^4_x  }\|  w\|_{L^2 L^4_x (\R_+)} + \| \na v\|_{L^2_tL^4_x  }\|  v\|_{L^2 L^4_x (\R_+)} \big) \\
& \leq c  \big( \| \na w\|_{ L^4  }\|  w\|_{  L^4 (\R_+)}  + \| \na w\|_{ L^4  }\|  v\|_{ L^4 (\R_+)}\\
  & \qquad + \| \na v\|_{ L^4  }\|  w\|_{  L^4 (\R_+)} + \| \na v\|_{ L^4  }\|  v\|_{  L^4 (\R_+)} \big).
\end{align*}
Since $r > 8$, from \eqref{0920-4}, we have $v \in L^\infty_t L^2_x$. Hence, $v$ is weak solution, i.e. $v \in  L^\infty_t L^2_x  \cap L^2_t H^1_x$.

We take $s_0$ with $\frac{n+2}2< s_0 < s$ such that $\frac1r <\frac1{r_1}:  = \frac1{s_0} -\frac1s$. We note that $ s_0 < r_1 < r$.
Then, we have
\begin{align}\label{1101-1}
\begin{split}
&\|\na^2  v\|_{L^{s_0}}  + \|D_t v\|_{L^{s_0}} + \|\na \pi \|_{L^{s_0}}\\
 \leq & c     \|(\na  w)w \|_{L^{s_0}} + \|(\na  w)v \|_{L^{s_0}} + \|(\na  v)w \|_{L^{s_0}} + \|(\na  v)v \|_{L^{s_0}}\\
 \leq & c     \|\na  w \|_{L^{s}} \|  w \|_{L^{r_1}} + \|\na  w \|_{L^{s}} \|  v \|_{L^{r_1}} + \|\na  v \|_{L^{s}} \|  w \|_{L^{r_1}} + \|\na  v \|_{L^{s}} \|  v \|_{L^{r_1}} < \infty.
\end{split}
\end{align}
We note that
\begin{align}\label{1101-2}
\begin{split}
w, v \in L^\infty (A \times (0,1)),\qquad \na v \in L^{s_0} (A \times (0,1)),\\
 \na w \in L^{r_-} (A \times (0,1))\quad\mbox{ for all }\,\,r_- < r_0.
\end{split}
\end{align}
Since $5<r_0$, we take $r_- $ with $5<r_- < r_0$ such that ${\rm div } F \in L^{r_-} (A \times (0,1))$. Then, applying \eqref{1101-1} and \eqref{1101-2} in the proof of  Theorem 1.5 in \cite{CK0206}, we have
\begin{align}
\na v \in C^{\frac12( 1 -\frac{5}{r_-})}_t C^{1 -\frac{5}{r_-}}_x(A \times (0,1)).
\end{align}
This implies that $\na v \in L^{r_0}(A \times (0,1)).$
We then set $u:=v+w$ and $p =\pi + q$, which becomes a weak solution of the
Navier-Stokes equations in $\R_+ \times (0, \infty)$. However,
\[
\| \na u\|_{L^{r_0}(A \times (0, 1))} \ge \| \na w\|_{L^{r_0}(A \times (0, 1))} -\| \na v\|_{L^{r_0}(A \times (0, 1))} \ge \| \na w\|_{L^{r_0}(A \times (0, 1))} -c.
\]
The righthand side becomes unbounded, and thus we obtain
$\| \na u\|_{L^{r_0}(A \times (0, 1))} = \infty$.
 \qed

\section{Appendix}
\setcounter{equation}{0}

\subsection{Proof of Lemma \ref{lemmaappendix}}
\label{alter-thm11}

Since $ \frac12 (x_n + y_n)^2 \leq x_n^2 + y_n^2 \leq 2 (x_n + y_n)^2$, for $t  > \frac12$, we note that
\begin{align}\label{0413-4}
\nonumber &\int_{\frac12}^t \int_0^2 y_n^{-\be} (\tau-\frac12)^{-\al}  (t -\tau)^\ga  e^{-\frac{(x_n + y_n)^2}{t-\tau}} dy_n d\tau\\
\nonumber  = & \int_{\frac12}^t  (\tau-\frac12)^{-\al}  (t -\tau)^{\frac12 -\frac{\be}2 +\ga} e^{-\frac{x_n^2}{t-\tau}} \int_0^{\frac2{\sqrt{t-\tau}}} y_n^{-\be} e^{-y_n^2} dy_n d\tau\\
\nonumber \approx  &\int_{\frac12}^t  (\tau-\frac12)^{-\al}  (t -\tau)^{\frac12 -\frac{\be}2 +\ga}  e^{-\frac{x_n^2}{t-\tau}} d\tau\\
 = &  \int_{\frac12}^{\frac12 ( t +\frac12)} \cdots d\tau+ \int_{\frac12 ( t +\frac12)}^t \cdots d\tau:=
 I_1 + I_2,
\end{align}
where we used that $ \int_0^{\frac1{\sqrt{t-\tau}}} y_n^{-\be} e^{-y_n^2} dy_n d\tau
 \approx 1$, since $0\le t-\tau\le 1$.
It is direct that $\frac12 ( t -\frac12) \leq t -\tau \leq  t -\frac12$ of $ \frac12 < \tau < \frac12 (t +\frac12 )$. Therefore, the term $I$ can be estimated as
\[
c (t -\frac12)^{\frac12 -\frac{\be}2 +\ga }  e^{-\frac{x_n^2}{t-\frac12 }}  \int_{\frac12}^{\frac12 ( t +\frac12)}  (\tau-\frac12)^{-\al}   d\tau   \le I_1  \le   c (t -\frac12)^{\frac12 -\frac{\be}2 +\ga }  e^{-\frac{x_n^2}{t-\frac12 }}  \int_{\frac12}^{\frac12 ( t +\frac12)}  (\tau-\frac12)^{-\al}   d\tau .
\]
Since $\int_{\frac12}^{\frac12 ( t +\frac12)}  (\tau-\frac12)^{-\al}   d\tau  \approx t^{1-\al}$,
 we obtain
\begin{align}\label{0620-1}
 c (t -\frac12)^{\frac32  -\al -\frac{\be}2 +\ga}  e^{-2\frac{x_n^2}{t-\frac12 }}\le I_1 \le c (t -\frac12)^{\frac32  -\al -\frac{\be}2 +\ga}  e^{-\frac{x_n^2}{t-\frac12 }}.
\end{align}
On the other hand, noting that $\frac12 ( t -\frac12) \leq  \tau -\frac12 \leq  t -\frac12$ of $  \frac12 (t +\frac12 ) < \tau < t$, we have
\begin{align*}
I_2 & \approx   c  (t-\frac12)^{-\al}\int_{\frac12 ( t +\frac12)}^t   (t -\tau )^{\frac12 -\frac{\be}2 +\ga}  e^{-\frac{x_n^2}{t-\tau}} d\tau\\
& =   c  (t-\frac12)^{-\al}\int_0^{\frac12 ( t -\frac12)}   \tau^{\frac12 -\frac{\be}2 +\ga}  e^{-\frac{x_n^2}{\tau}} d\tau \\
&=  c  (t-\frac12)^{-\al} x_n^{3 -\be +2\ga} \int^\infty_{{\frac{2 x_n^2}{t-\frac12}}  }   s^{-\frac52 +\frac{\be}2-\ga }  e^{-\tau} d\tau,
\end{align*}
where we used
the change of variables ( $r = \frac{x_n^2}{s}$) in the last equality.

{$\bullet$ \,\,( Case $ \ga  -\frac{\be}2 \ge  -\frac32 $)}\,\,We treat the both cases of
$ \frac{2 x_n^2}{t-\frac12} \leq 1$ and $ \frac{2 x_n^2}{t-\frac12} \ge 1$ separately.
Firstly, if  $ \frac{2 x_n^2}{t-\frac12} \leq 1$, it follows that
\[
I_2  = c  (t-\frac12)^{-\al} x_n^{3 -\be +2\ga} \int^1_{{\frac{2 x_n^2}{t-\frac12}}  }   s^{-\frac52 +\frac{\be}2 -\ga}  e^{-\tau} d\tau + c  (t-\frac12)^{-\al} x_n^{3 -\be +2\ga} \int_1^\infty s^{-\frac52 +\frac{\be}2 -\ga}  e^{-\tau} d\tau
\]
\[
\approx   c  (t-\frac12)^{-\al} x_n^{3 -\be +2\ga} \big( {\frac{2 x_n^2}{t-\frac12}}  \big)^{-\frac32 + \frac{\be}2 -\ga}  + c  (t-\frac12)^{-\al} x_n^{3 -\be +2\ga}
\]
\begin{equation}\label{0620-2}
\approx   c  (t-\frac12)^{\frac32 -\frac\be{2}-\al +\ga} \approx   c  (t-\frac12)^{\frac32 -\frac\be{2}-\al +\ga} e^{-\frac{2 x_n^2}{t-\frac12} },
\end{equation}
where we used $e^{-1}\le e^{-\frac{x_n^2}{t-\frac12}}\le 1$.

Next, we consider the case $ \frac{2 x_n^2}{t-\frac12} \geq 1$.
We note via the l'Hospital's Theorem that
\begin{align}\label{lhospital}
\lim_{a \ri \infty} \frac{\int_a^\infty \tau^{-\frac52 +\frac{\be}2 -\ga} e^{-\tau}d\tau }{ a^{-\frac52 +\frac{\be}2-\ga} e^{-a} } = \lim_{a \ri \infty} \frac{-a^{-\frac52 +\frac{\be}2 -\ga} e^{-a} }{ ( -\frac52 +\frac{\be}2-\ga) e^{-a} - a^{-\frac52 +\frac{\be}2 -\ga} e^{-a} } =1.
\end{align}
Hence, due to \eqref{lhospital}, it is straightforward that
\[
I_2 \approx  (t-\frac12)^{\frac{3}{2}-\al-\frac{\beta}{2}+\gamma}   e^{-{\frac{2 x_n^2}{t-\frac12}}}.
\]

$\bullet$\,\,
{\bf (Case $  \ga -\frac{\be}2  <  -\frac32 $)}\,\,
In case that $ \frac{2 x_n^2}{t-\frac12} \leq 1$, it is direct that
\begin{align*}
I_2
& \approx   c  (t-\frac12)^{-\al} x_n^{3 -\be +2\ga}.
\end{align*}
If $ \frac{2 x_n^2}{t-\frac12} \geq 1$, it follows from \eqref{lhospital} that
\begin{align}\label{0620-3}
I_2 & \approx  c  (t-\frac12)^{-\al} x_n^{3 -\be +2\ga} ({\frac{ x_n^2}{t-\frac12}}  )^{-\frac52 +\frac{\be}2 -\ga} e^{-{\frac{2 x_n^2}{t-\frac12}} } \\
\notag& \approx   c  (t-\frac12)^{-\al} x_n^{3 -\be +2\ga}  e^{-{\frac{2 x_n^2}{t-\frac12}}  }.
\end{align}
 Hence, from \eqref{0620-1}, \eqref{0620-2}, \eqref{0620-3} and \eqref{0620-3}, we  complete the proof of  Lemma  \eqref{lemmaappendix}.

\subsection{Proof of Lemma \ref{lemma0709-1}}
\label{prooflemma0709-1}

First, we assume that $l =0$ and $D^k_{x'} = D^{k-1}_{x'} D_{x_i}, \,\, k \geq 1, \,\, 1 \leq i \leq n-1 $. So,
\begin{align*}
 D^k_{x'} \int_{\Rn}    \Ga'(x'-z',t)    N( z',x_n) dz' =   D^{k-1}_{x'} \int_{\Rn}    \Ga'(x'-z',t)   D_{z_i} N( z',x_n) dz'.
\end{align*}

We divide $\Rn$ by three disjoint sets
$D_1, D_2$ and $D_3$, which are defined by
\[
D_1=\bket{z'\in\Rn: |x'-z'| \leq \frac1{10} |x'|},
\]
\[
D_2=\bket{z'\in\Rn: |z'| \leq \frac1{10} |x'|}, \qquad D_3=\Rn\setminus
(D_1\cup D_2).
\]
We then split the following integral into three terms as follows:
\begin{equation}\label{0730-2}
\int_{\Rn}D^{k-1}_{z'} \Ga'(x'-z',t)   D_{z_i}  N( z', x_n) dz' =
\int_{D_1}\cdots + \int_{D_2} \cdots+ \int_{D_3}\cdots := J_1 + J_2
+J_3.
\end{equation}

Noting that $ |D^{k}_{z'} \Ga'(x'-z',t) | \leq ct^{-\frac{n+k-1}2 } e^{-\frac{|x'|^2}{2t}}$  for $z' \in D_2$ and $\int_{D_2} D_{z_i}  N( z',x_n) dz' =0$,  we have
\begin{align}\label{estJ2}
\notag |J_2| & =\abs{\int_{D_2} D_{z_i}  N( z',x_n)  \big(   D^k_{z'}\Ga'(x'-z',t)  - D^{k}_{x'}\Ga'(x',t) \big) dz'
}\\
&\notag \leq   ct^{-\frac{n -1}2 -\frac{k}2   } e^{-\frac{|x'|^2}{2t}} |    \int_{D_2}
 \frac{1}{  |z'|^{n-2} }  dz'\\
& \leq   ct^{-\frac{n+k-1}2  }|x'|  e^{-\frac{|x'|^2}{2t}}\le c|x'|^{-\frac{n+k-1}2  }t^{\frac{1}{2}}\le ct^{\frac{1}{2}},
\end{align}
where the mean value theorem is used. Here we also used that $e^{-\frac{|x'|^2}{2t}}\le c(\frac{|x'|^2}{t})^{-\frac{n+k}{2}}$.
On the other hand, the term $J_3$ is controlled as follows:
\begin{align}\label{est-J3}
 |J_3|
  \leq   \frac{c}{|x'|^{n-1}  }  t^{-\frac{n-1}2 -\frac{k-1}2}  \int_{\{|z'| \geq \frac1{10} |x'|\}}  e^{-\frac{|z'|^2}{t}} dz'
 \leq  \frac{c}{|x'|^{n-1}  }  t^{-\frac{k-1}2} e^{-\frac{|x'|^2}t}\le ct^{\frac{1}{2}},
\end{align}
where we used that $e^{-\frac{|x'|^2}{2t}}\le c(\frac{|x'|^2}{t})^{-\frac{k}{2}}$.

Now, it remains to estimate $J_1$.
Due to the integration by parts, it follows that
\begin{align*}
 J_{1} & = \int_{D_1} D^{k-1}_{x'} \Ga'(x' -z',t)  D_{z_i} N(z',x_n)  dz'\\
 & = (-1)^{k-1} \int_{D_1} D^{k-1}_{z'} \Ga'(x' -z',t)  D_{z_i} N(z',x_n)  dz'\\
 & = (-1)^{k-1} \sum_{1 \leq k' \leq k-2} (-1)^{k'} \int_{\pa D_1}   D^{k'}_{z'} \Ga'(x' -z',t)  D^{k-1 -k'}_{z'}D_{z_i} N(z',x_n) {\bf n}_{k'} d\sigma(z')\\
&  + \int_{D_1}  \Ga'(x' -z',t)  D^{k-1}_{z'} D_{z_i} N(z',x_n)  dz'\\
& := J_{11} + J_{12}.
 \end{align*}
Since the magnitude of $z'\in \partial D_1$ is comparable to $|x'|$, $J_{11}$ is controlled as follows:
\begin{align}\label{est-J11}
\notag |J_{11}|
& \leq c \sum_{1 \leq k' \leq k-2}     t^{-\frac{n-1}2 -\frac{k'}2}e^{-\frac{|x'|^2}t} \int_{\pa D_1}  |  D^{k-1 -k'}_{z'}D_{z_i} N(z',x_n) {\bf n}_{k'}| d\sigma(z')\\
& \notag \leq c \sum_{1 \leq k' \leq k-2}     t^{-\frac{n-1}2 -\frac{k'}2}e^{-\frac{|x'|^2}t} \int_{\pa D_1} |x'|^{-n  -k + k'}  d\sigma(z')\\
& \leq c \sum_{1 \leq k' \leq k-2}     t^{-\frac{n-1}2 -\frac{k'}2} |x'|^{ -k +1 + k'}   e^{-\frac{|x'|^2}t}\le ct^{\frac{1}{2}},
\end{align}
where we used that $e^{-\frac{|x'|^2}{t}}\le c(\frac{|x'|^2}{t})^{-\frac{n+k'}{2}}$.
Meanwhile, we decompose $J_{12}$ in the following way.
 \begin{align*}
\notag  J_{12} & =     \int_{ D_1}    \Ga'(x' -z',t)  D^{k-1 }_{z'}D_{z_i} N(z',x_n) dz'
\\
\notag &
 =    \int_{\{|z'| \leq    \frac1{10}\frac{ |x'|}{\sqrt{t}}\}}     \Ga'(z',1) \Big( D^{k-1}_{z'} D_{x_i} N(x' -\sqrt{t}z',x_n) - D^{k-1}_{x'}D_{x_i} N(x',x_n) \Big)dz'\\
\notag &\quad + D^{k-1}_{x'} D_{x_i} N(x',x_n) \int_{\Rn}  \Ga'(z',1)  dz' - D^{k-1}_{x'} D_{x_i} N(x',x_n) \int_{\{|z'| \geq  \frac1{10}\frac{ |x'|}{\sqrt{t}}\}}  \Ga'(z',1)  dz'\\
& = J_{121} +D^{k-1}_{x'} D_{x_i} N(x',x_n) + J_{122},
\end{align*}
where we used $\int_{\Rn}  \Ga'(z',1)  dz'=1$.
Since $\frac{\sqrt{t}}{|x'|}\le 1$, we observe that
\begin{equation}\label{0730-1}
|J_{121} (x',t)| \leq c |x'|^{-n -k+1} t^\frac12   \int_{\{|z'| \leq  \frac1{10}\frac{ |x'|}{\sqrt{t}}\}}    e^{-|z'|^2}|z'|   dz'\le ct^{\frac{1}{2}  },
\end{equation}
\begin{equation}\label{est-J122}
|J_{122} (x',t)| \leq c |x'|^{-n -k +2}    \int_{\{|z'| \geq  \frac1{10}\frac{ |x'|}{\sqrt{t}}\}}    e^{-|z'|^2}   dz'\le c e^{-\frac{|x'|^2}{2t}}\le ct^{\frac{1}{2}  }.
\end{equation}
Setting $J= J_2 + J_3 + J_{11} + J_{121} + J_{122}$ and adding up \eqref{estJ2} -\eqref{est-J122}, we deduce   \eqref{0515-1} for $k \geq 1$ and $l =0$.

Next, we consider the case that normal derivative is taken into account, i.e. $l\ge 1$.
We note first that $N(z', x_n)$ is regular in the regions of $D_1$ and $D_3$, where all order of normal derivatives can be directly applied to $N(z', x_n)$. Therefore, $J_1$ and $J_3$ in \eqref{0730-2} can be computed similarly as in the above case that $l=0$.
Since its verifications are just tedious  repetitions, it suffices that we estimate only $J_2$ in \eqref{0730-2}.

Firstly, in case that $l =1$,
\begin{align}\label{est-J2}
 |J_2|  \leq  c t^{-\frac{n-1}2-\frac{k}2} e^{-\frac{|x'|^2}t}\int_{|z'| \leq \frac1{10} |x'|} \frac{x_n}{(|z'|^2 + x_n^2)^{\frac{n}2}}  dz'
 \leq ct^{\frac{1}{2}}.
\end{align}
Secondly,
if $l =2$, then it follows due to $D_{x_n}^2 N(z', x_n)=-\Delta'_{z'}N(z', x_n)$ that
\begin{align}\label{est-J2-1}
\notag J_2 & = c \int_{|z'| \leq \frac1{10} |x'|} D_{x'}^k \Ga'(x'-z', t) D_{x_n} D_{x_n} N(z',x_n)  dz'\\
& \notag = - c \int_{|z'| \leq \frac1{10} |x'|} D_{x'}^k \Ga'(x'-z', t) \De'_{z'} N(z',x_n)  dz'\\
& \notag =   c \int_{|z'| \leq \frac1{10} |x'|} \na_{z'} D_{x'}^k \Ga'(x'- z', t) \cdot  \na_{z'} N(z',x_n)  dz' \\
&\notag\quad -  c \int_{|z'| =  \frac1{10} |x'|}  D_{x'}^k \Ga'(x'-z', t)  \na_{z'} N(z',x_n) \cdot \nu' dz'\\
& = J_{1a} + J_{1b}.
\end{align}
Recalling that $\int_{|z'| \leq \frac1{10} |x'|}  \na_{z'} N(z',x_n)  dz'=0$, we observe that
\begin{align*}
|J_{1a}| & = c \abs{\int_{|z'| \leq \frac1{10} |x'|} \big( \na_{z'} D_{x'}^k \Ga'(x'- z', t) - \na_{x'} D_{x'}^k  \Ga(x', t) \big) \cdot  \na_{z'} N(z',x_n)  dz'}\\
&  \leq   ct^{-\frac{n -1}2 -\frac{k+2}2   } e^{-\frac{|x'|^2}{2t}} \int_{D_2}
 \frac{1}{  |z'|^{n-2} }  dz'\\
& \leq   ct^{-\frac{n-1}2 -\frac{k+2}2  }|x'|  e^{-\frac{|x'|^2}{2t}}
\le ct^{\frac{1}2   },
\end{align*}
where $  e^{-\frac{|x'|^2}{2t}}\le c(\frac{|x'|^2}{t})^{-\frac{n+k-2}{2}}$ is used.
The second term can be similarly estimated as follows:
\begin{align*}
|J_{2b} | \leq ct^{-\frac{n -1}2 -\frac{k}2   } e^{-\frac{|x'|^2}{2t}} |x'|^{-1}
\le ct^{\frac{1}2   }.
\end{align*}
For the case $l>2$, we can convert the order of normal derivatives to tangential derivatives by using $D_{x_n}^{2m} N(z', x_n)=(-1)^m(\Delta'_{z'})^mN(z', x_n)$ for $m\le [\frac{l}{2}]$, which reduces it to the case either $l=0$ or $l=1$.
Repeating the above processing similarly, we can have  \eqref{0515-1}-\eqref{Jkl}. Since its computations are rather straightforward, we skip its details. We deduce the lemma.

\subsection{A Figure of disjoint sets}\label{boundaryset-1}

The two dimensional figure of disjoint sets $A_{i1}, A_{i2}, B_{i1}$ and $B_{i2}$ are given as follows:

\begin{tikzpicture}[scale=3]

\draw (-1, -1) node {$A_{11}$};
\draw (1, 1) node {$A_{11}$};

\draw (1, -1) node {$A_{12}$};
\draw (-1, 1) node {$A_{12}$};

\draw (2.05,-0.099) node {$x_1$};
\draw(0.13,1.5) node {$x_2$};

\draw (-2,0) -- (2,0);
\draw (0,-1.6) -- (0, 1.6);

\draw (0, 1.6)  -- (-0.05, 1.55)    ;
\draw (0, 1.6)  -- (0.05, 1.55);

\draw (2,0) -- (1.95,0.05);
\draw (2,0) -- (1.95,-0.05);

\draw (-1.7, 0.3) -- (-0.62, 0.11 );
\draw (-1.7, -0.3) -- (-0.62, -0.11 );

\draw (1.7, 0.3) -- (0.62, 0.11 );
\draw (1.7, -0.3) -- (0.62, -0.11 );

\draw (0.3,-1.5) -- (0.127, -0.63);
\draw (-0.3,-1.5) -- (-0.127, -0.63);

\draw (0.3,1.5) -- (0.127, 0.63);
\draw (-0.3,1.5) -- (-0.127, 0.63);

\draw (-1.0,1.25) -- (-0.5,0.63);
\draw (-1.4,0.9) -- (-0.63,0.5);

\draw (0.5,-0.63) -- (1.0,-1.25);
\draw (0.63,-0.5) -- (1.4,-0.9);

\draw (-0.5,-0.63) -- (-1.0,-1.25);
\draw (-0.63,-0.5) -- (-1.4,-0.9);

\draw (1.0,1.25) -- (0.5,0.63);
\draw (1.4,0.9) -- (0.63,0.5);

\draw (0.63, 0.63) -- (0.63, -0.63);
\draw (-0.63, -0.63) -- (0.63, -0.63);
\draw (-0.63, -0.63) -- (-0.63, 0.63);
\draw (-0.63, 0.63) -- (0.63, 0.63);

\draw (-0.4, -0.1) node {$(-2, 0)$};
\draw (0.45,0.1) node {$(2, 0)$};
\draw (-0.2, -0.5) node {$(0, -2)$};
\draw (0.15, 0.5) node {$(0, 2)$};

\draw (0.1, -0.1) node {$ 0$};

\draw (-1.5, 0) node {$B_{11}$};
\draw (1.5, 0) node {$B_{11}$};

\draw (0, -1.3) node {$B_{12}$};
\draw (0, 1.3) node {$B_{12}$};

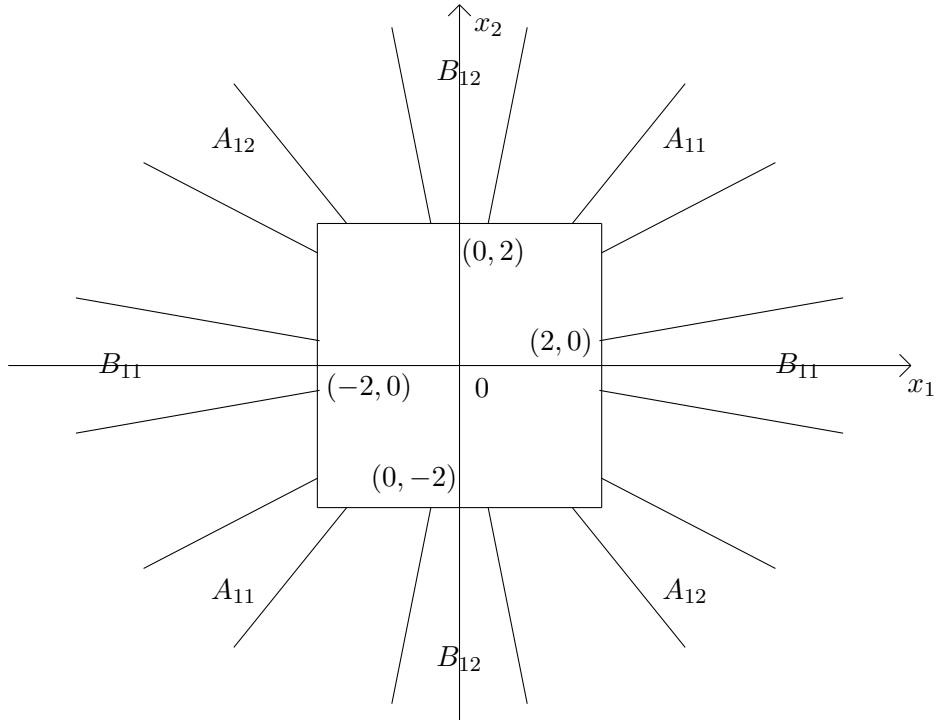
\draw (0.0, -1.7) node {Figure; $A_{1j}$ and $B_{1j}$ in ${\mathbb R}^2$, $ j =1, 2$ };

\end{tikzpicture}

\section*{Acknowledgement}
T. Chang is partially supported by NRF-2020R1A2C1A01102531. K. Kang is supported
by NRF-2019R1A2C1084685.

\begin{equation*}
\left.
\begin{array}{cc}
{\mbox{Tongkeun Chang}}\qquad&\qquad {\mbox{Kyungkeun Kang}}\\
{\mbox{Department of Mathematics }}\qquad&\qquad
 {\mbox{Department of Mathematics}} \\
{\mbox{Yonsei University
}}\qquad&\qquad{\mbox{Yonsei University}}\\
{\mbox{Seoul, Republic of Korea}}\qquad&\qquad{\mbox{Seoul, Republic of Korea}}\\
{\mbox{chang7357@yonsei.ac.kr }}\qquad&\qquad
{\mbox{kkang@yonsei.ac.kr }}
\end{array}\right.
\end{equation*}


\begin{thebibliography}{10}

%\bibitem{amman-anisotropic} H. Amann, {\it
%Anisotropic function spaces and maximal regularity for parabolic problems. Part 1.
%Function spaces,}  Jindřich Nečas Center for Mathematical Modeling Lecture Notes, 6. Matfyzpress, Prague, vi+141(2009).

\bibitem{BL}  J. Bergh and J. L\"ofstr\"om, {\it Interpolation Spaces, An Introduction},
Springer-Verlag, Berlin 1976.

\bibitem{CJ} T. Chang and B.  Jin, {\it Initial and boundary values for $L^q_\al(L p )$  solution
of the Navier-Stokes equations in the half-space}, J. Math. Anal. Appl.  {\bf 439},  no. 1, 70-90 (2016).

%\bibitem{CJ2} T. Chang and B. Jin, {\it Global in time solvability of the Navier-Stokes equations in the half-space }. J. Differential Equations, {\bf  267}, no. 7, 4293–4319  (2019).

%\bibitem{CK} T. Chang and K. Kang, {\it Estimates of anisotropic Sobolev spaces with mixed norms for the Stokes system in a half-space}, Ann. Univ. Ferrara Sez. VII Sci. Mat. {\bf  64},  no. 1, 47-82   (2018).

%\bibitem{CKsima} T. Chang, H. Choe and K. Kang, {\it On maximum modulus estimates of the Navier-Stokes equations with non-zero boundary data},  SIAM J. Math. Anal. {\bf 50}, no. 3, 3147-3171(2018).

\bibitem{CKca} T. Chang and K. Kang, {\it On Caccioppoli's inequalities of Stokes equations and Navier-Stokes equations near boundary},  J. Differential Equations, {\bf 269}, no. 9, 6732–6757 (2020).

\bibitem{CK0206} T. Chang and K. Kang, {\it Local regularity near boundary for the Stokes anad Navier-Stokes equations}, Preprint: arXiv:2110.07162.


%\bibitem{CSYT08} C.-C. Chen, R. M. Strain, H.-T. Yau and T.-P. Tsai,
%{\it Lower bound on the blow-up rate of the axisymmetric
%Navier-Stokes equations}, Int. Math. Res. Not. IMRN  {\bf 2008}, no.
%9, pp 31.

%\bibitem{Choe-Yang15} H. Choe and M. Yang, {\it Hausdorff measure of the singular set in
%the incompressible magnetohydrodynamic equations}, Comm. Math. Phys.
%{\bf 336}, no. 1, 171-198(2015).

%\bibitem{CM2006} F. Crispo and P. Maremonti, {\it On the $(x,t)$ asymptotic properties of solutions of the Navier-Stokes equations in the half-space}. Zap. Nauchn. Sem. S.-Peterburg. Otdel. Mat. Inst. Steklov. (POMI) 318 (2004), Kraev. Zadachi Mat. Fiz. i Smezh. Vopr. Teor. Funkts. 36 [35], 147–202, 311; reprinted in J. Math. Sci. (N.Y.) 136 (2006), no. 2, 3735–3767.

%\bibitem{DT} H. Dappa and  H. Triebel, {\it On anisotropic Besov and
%                  Bessel Potential spaces}, Approximation and function spaces, 69-87 (Warsaw, 1986), Banach Center
%                   Publ., 22, PWN, Warsaw(1989).


%\bibitem{galdi} G.P. Galdi, {\it An Introduction to the Mathematical Theory of the Navier-Stokes Equations}, vol. I, linearlised steady problems, Springer Tracts in Natural Philosopy vol. 38, Springer 1994.


%\bibitem{Giaquinta-Modica} M. Giaquinta and G. Modica, {\it Nonlinear systems of the type of the stationary Navier-Stokes
%system}, J. Reine Angew. Math.  {\bf 330},  173-214 (1982).

%\bibitem{GGS} M. Giga, Y. Giga and H. Sohr, {\it $L^p$ estimates for the Stokes system},
%Lect. Notes Math. {\bf 1540},  (1993), 55-67.

\bibitem{JK} D. Jerison and C. E. Kenig, {\it The inhomogeneous Dirichlet Problem in Lipschitz domains},
J. of Funct. Anal, 130, 161-219(1995).


%\bibitem{J}  B.  Jin, {\it On the Caccioppoli inequality of the unsteady Stokes system},  Int. J. Numer. Anal. Model. Ser. B, {\bf  4},  no. 3, 215-223(2013).

%\bibitem{Jin-Kang17}  B.  Jin and K. Kang, {\it Caccioppoli type inequality for non-Newtonian Stokes system and a
%local energy inequality of non-Newtonian Navier-Stokes equations
%without pressure}, Discrete Contin. Dyn. Syst. {\bf  37}, no. 9,
%4815-4834(2017).

\bibitem{Kang05} K. Kang, {\it Unbounded normal derivative for the Stokes system
near boundary}, Math. Ann.  {\bf 331} no. 1, 87-109(2005).

%\bibitem{Kang04} K. Kang, {\it On regularity of stationary Stokes and Navier-Stokes
%equations near boundary}, J. Math. Fluid Mech.  {\bf 6} no. 1,
%78-101(2004).

%\bibitem{KLLT21-G} K. Kang, B. Lai, C.-C. Lai, T.-P. Tsai, {\it The Green tensor of the nonstationary
%Stokes system in the half space}, Preprint:  https://arxiv.org/abs/2011.00134.

\bibitem{KLLT21} K. Kang, B. Lai, C.-C. Lai, T.-P. Tsai, {\it Finite energy Navier-Stokes flows with unbounded gradients induced by localized flux in the half-space}, accepted in  Trans. Amer. Math. Soc.

\bibitem{KS} H. Koch and V.A. Solonnikov, {\it $L_p$estimates for a solution to the nonstationary Stokes equations. Function theory and phase transitions}, J. Math. Sci. (New York)  {\bf 106},    no. 3, 3042-3072 (2001).

\bibitem{KS2002} H. Koch and V.A. Solonnikov, {\it Lq-estimates of the first-order derivatives of solutions to the nonstationary Stokes problem}, Nonlinear problems in mathematical physics and related topics, I, 203–218, Int. Math. Ser. (N. Y.), 1, Kluwer/Plenum, New York, 2002.

%\bibitem{LS}
%O.A. Lady$\check{\rm z}$enskaja, V.A. Solonnikov and N.N. Ura${\rm l}$ceva, {\it Linear and Quasilinear Equations of Parabolic Type}, (Russian) Translated from the Russian by S. Smith. Translations of Mathematical Monographs, Vol. 23 American Mathematical Society, Providence, R.I. 1968.


%\bibitem{Seregin00} G. A. Seregin, {\it Some estimates near the boundary for solutions to the non-stationary
%linearized Navier-Stokes equations}, (English, Russian summary) Zap.
%Nauchn. Sem. S.-Peterburg. Otdel. Mat. Inst. Steklov. (POMI)  {\bf
%271} (2000),  Kraev. Zadachi Mat. Fiz. i Smezh. Vopr. Teor. Funkts.
%31, 204--223, 317;  translation in J. Math. Sci. (N. Y.)  {\bf 115}
%no. 6,  2820-2831(2003).

\bibitem{Seregin-Sverak10} G. A. Seregin and V. \u{S}ver\'ak, {\it On a bounded shear flow in half-space}, Zap. Nauchn. Sem. S.-Peterburg. Otdel. Mat. Inst. Steklov. (POMI) 385 (2010), Kraevye Zadachi Matematicheskoĭ Fiziki i Smezhnye Voprosy Teorii Funktsiĭ. 41, 200–205, 236; J. Math. Sci. (New York)  {\bf 178},   no. 3, 353–356 (2011).

%\bibitem{So}  V. A. Solonnikov, {\it Estimates of the solutions of the nonstationary Navier-Stokes system.
%(Russian) Boundary value problems of mathematical physics and related questions in the theory of functions},
%$7$. Zap. Naucn. Sem. Leningrad. Otdel. Mat. Inst. Steklov. (LoMI) 38: 153-231. Translated in J. Soviet Math. 1977, 8: 47-529.

\bibitem{So1} V. A. Solonnikov, {\it On estimates of the solutions of the non-stationary Stokes problem in anisotroptc Sobolev spaces and on estimates for the resolvent of the Stokes operator}. Russian Math. Surveys 58:2 331-365(2003).



%\bibitem{S02} V. A. Solonnikov, {\it The initial boundary-value problem for a generalized Stokes system in a half-space}, Journal of Mathematical Sciences, {\bf 115}, no. 6, 2832-2861(2003).

\bibitem{S03} V. A. Solonnikov, {\it Estimates for solutions of the nonstationary Stokes problem in anisotropic Sobolev spaces and estimates for the resolvent of the Stokes operator},
 (Russian) Uspekhi Mat. Nauk 58 (2003), no. 2(350), 123–156; translation in Russian Math. Surveys 58 (2003), no. 2, 331–365

%\bibitem{St} E.M. Stein, {\rm  Singular Integrals and Differentiability Properties of Functions,}
%Princeton Mathematical Series, No. 30 Princeton University Press, Princeton, N.J. 1970.


%\bibitem{Sverak-Tsai98} V. \v{S}ver\'ak and T.-P. Tsai, {\it On the spatial decay of 3-D
%steady-state Navier-Stokes flows},  Comm. Partial Differential
%Equations  {\bf 25} no. 11-12, 2107-2117(2000).

%\bibitem{Triebel2} H. Triebel, {\it  Theory of Function Spaces. III,} Monographs in Mathematics, 100. Birkh$\ddot{\rm a}$user Verlag, Basel, 2006.
%

%\bibitem{Tr} H. Triebel, {\it Interpolation Theory, Function Spaces, Differential Operators}, North-Holland Publishing company, 1978.

%\bibitem{Triebel} H. Triebel, {\it Theory of Function Spaces,} Monographs in Mathematics, 78. Birkh$\ddot{\rm a}$user Verlag, Basel, 1983.

%\bibitem{Wolf15} J. Wolf, {\it On the local regularity of suitable weak solutions to the
%generalized Navier-Stokes equations}, Ann. Univ. Ferrara Sez. VII
%Sci. Mat.  {\bf 61} no. 1, 149-171(2015).










\end{thebibliography}
\end{document}